\documentclass[]{amsart}
\usepackage{graphicx} % Required for inserting images
\usepackage{amsfonts,latexsym,rawfonts,amsmath,amssymb, amsthm,bm}
\usepackage{latexsym,lscape,rawfonts}
\usepackage[margin=3cm]{geometry}
\usepackage[backref=page]{hyperref}
\usepackage{cleveref}
\usepackage{orcidlink}
\usepackage[T1]{fontenc}

\usepackage{latexsym,amsmath,xcolor,multicol,booktabs,calligra}
\usepackage{amsfonts,latexsym,rawfonts,amsmath,amssymb, amsthm,bm}
\usepackage{graphicx,subcaption,pstricks,listings,stackengine}
\usepackage{float}
\usepackage{comment}
\usepackage{tikz-cd}
\usepackage{tikz}
\usetikzlibrary{arrows.meta,decorations.markings,calc}
\usepackage{circuitikz}

\newtheorem{thm}{Theorem}[section]
\newtheorem{cor}[thm]{Corollary}
\newtheorem{lem}[thm]{Lemma}
\newtheorem{prop}[thm]{Proposition}

\theoremstyle{definition}
\newtheorem{defn}[thm]{Definition}

\newtheorem{remark}[thm]{Remark}
\numberwithin{equation}{section}

\newcommand{\Z}{\mathbb{Z}}
\newcommand{\C}{\mathcal{C}}
\newcommand{\N}{\mathcal{N}}
\newcommand{\A}{\mathcal{A}}
\newcommand{\G}{\mathcal{G}}

\newcommand{\s}{\mathcal{S}}

\newcommand{\E}{\mathrm{E}}
\newcommand{\Ecal}{\mathcal{E}}
\newcommand{\M}{\mathrm{M}}
\newcommand{\V}{\mathrm{V}}
\newcommand{\Vcal}{\mathcal{V}}

\newcommand{\Cf}{\mathrm{CF}}
\newcommand{\CF}{\mathrm{CF}}
\newcommand{\Out}{\mathrm{Out}}

\newcommand{\Aut}{\mathrm{Aut}}
\newcommand{\SAut}{\mathrm{SAut}}
\newcommand{\Mod}{\mathrm{Mod}}
\newcommand{\Fact}{\mathrm{Fact}}
\newcommand{\Perm}{\mathrm{Perm}}
\newcommand{\FR}{\mathrm{FR}}
\newcommand{\PC}{\mathcal{PC}}
\newcommand{\id}{\mathrm{Id}}
\newcommand{\Isom}{\mathrm{Isom}}
\newcommand{\Fix}{\mathrm{Fix}}
\newcommand{\Fixdim}{\mathrm{Fixdim}}
\newcommand{\Niel}{\mathrm{Niel}}
\newcommand{\Post}{\mathrm{Post}}

\newcommand{\CAT}[1]{\mathrm{CAT(#1)}}
\newcommand{\generate}[1]{\langle #1 \rangle}
\newcommand{\Elp}{\mathrm{Elp}}

\title{Fixed points for actions of automorphism groups of free products on $\CAT0$ spaces}

\date{\today}
\author[Qian Chen]{Qian Chen\orcidlink{0009-0004-9630-9228}}	
\address{School of Mathematics and Statistics, Xi'an Jiaotong University, Xi'an 710049, P. R. China}	
\email{1161572435@stu.xjtu.edu.cn}

\author[Qiang Zhang]{Qiang Zhang\orcidlink{0000-0001-6332-5476}}	
\address{School of Mathematics and Statistics, Xi'an Jiaotong University, Xi'an 710049, P. R. China}	
\email{zhangq.math@mail.xjtu.edu.cn}

\date{\today}\author[Dongxiao Zhao]{Dongxiao Zhao\orcidlink{0009-0006-7919-9244}}	
\address{School of Mathematics and Statistics, Xi'an Jiaotong University, Xi'an 710049, P. R. China}	
\email{zdxmath@stu.xjtu.edu.cn}

\date{\today}
\thanks{The authors are partially supported by NSFC (No. 12471066).}

\begin{document}

\begin{abstract}
We study fixed-point properties of automorphism groups of free products acting on complete $\CAT0$ spaces. Under the assumption of certain fixed-point properties of the free factors, we obtain lower bounds on the smallest dimension in which automorphism groups can act on a complete $\CAT0$ space without a global fixed point. 
\end{abstract}

\maketitle

\section{Introduction}

Fixed-point properties of isometric group actions on non-positively curved spaces play an important role in geometric group theory. Given a metric space $X$, for any subgroup $G$ of the isometric group $\Isom(X)$, its fixed point set in $X$ is $$\Fix(G) := \{x \in X \mid g (x)=x,~\forall g \in G\}=\cap_{g\in G}\Fix(g).$$
We say that the \emph{action has a (global) fixed point} if $\Fix(G)\neq\emptyset$.

The classical fixed-point property for actions on trees is Serre's property $\rm FA$: a group has \emph{property $\rm FA$} if every isometric action on a simplicial tree has a fixed point. Culler and Vogtmann introduced the analogous \emph{property $\rm F\mathbb R$}, i.e., every isometric action on an $\mathbb R$-tree has a fixed point. Important examples of groups with such fixed-point properties include $\mathrm{SL}(n,\Z)$, $\Aut(F_n)$, $\Out(F_n)$, for $n\geq 3$. See \cite{Bogopolski,Culler,Serre}.

Farb \cite{Farb(Helly'stheorem)} extended property $\rm FA$ to actions on higher-dimensional non-positively curved spaces. A group $G$ is said to have \emph{property $\mathrm{FA}_n$} (resp. \emph{strong $\mathrm{FA}_n$ property}) if every isometric (resp. semi-simple isometric) action of $G$ on any $n$-dimensional, $\CAT0$ cell complex (resp. complete $\CAT0$ space) has a fixed point. In particular, simplicial trees and complete $\mathbb R$-trees are $\CAT0$ spaces of dimension at most one.

Throughout this paper, we adopt the following notion of topological dimension, as used by Bridson:
\begin{defn}\cite[Definition 2.1]{Bridson(Hellytype)}
    A topological space $X$ has \emph{dimension at most $d$}, written $\dim(X)\leq d$, if for every closed subspace $K\subseteq X$, every continuous map $f:K\to \mathbb S^r$ to a sphere of dimension $r\geq d$ extends to a continuous map $X\to \mathbb S^r$.
\end{defn}

Using the topology of nerves of convex subsets together with a bootstrap technique, Bridson \cite{Bridson(semisimple2), Bridson(MCG)} proved that the mapping class group of an orientable closed surface of genus larger than $n$ has strong $\mathrm{FA}_n$ property. Similar discussion can be found in matrix groups \cite{Ye1,Ye2}, Thompson groups \cite{Ka1,Ka2} and Coxeter groups \cite{Barnhill}. 

In recent years, several automorphism groups of free products have been studied: Leder \cite{Leder} and Andrew \cite{Andrew} studied actions of automorphism groups of free products on trees. They showed that when sufficiently many isomorphic factors occur in the free product decomposition, the corresponding automorphism group has property $\rm{FA}$.
From another perspective, Varghese \cite{OV(Aut(F_n))} studied the isometric actions of $\Aut(F_n)$ on $\CAT0$ spaces; Bridson \cite{Bridson(Hellytype)} showed the following result. 
\begin{thm}\cite[Theorem A]{Bridson(Hellytype)}\label{Bridson theorem}
    If $n\geq 3m$ and $d<2m$, or $n\geq 3m+2$ and $d<2m+1$, then $\Aut(F_n)$ has a fixed point whenever it acts by isometries on a complete $\CAT0$ space of dimension $d$.
\end{thm}

In this paper, we investigate the automorphism groups of free products acting by isometries on complete $\CAT0$ spaces (with dimension $\geq 1$). Since infinitely generated countable groups do not have Serre's property $\rm FA$, we restrict our attention to finitely generated free products.
Relying on two key properties (Properties (I) and (II), see Section \ref{sect 3}), we have the following theorem. 

\begin{thm}\label{main thm}
    Let $G=\ast_{i=1}^{m}\ast_{j=1}^{\ell_i} G_{i,j}$ be a free product as in Eq. (\ref{free decomposition}) and let $\Aut(G)$ act by isometries on a complete $\CAT0$ space $X$ with $\dim(X)\leq d$. If the following hold:
    \begin{itemize}
        \item each $\ell_i\geq 2d+2$;
        \item $G$ satisfies Properties (I) and (II),
    \end{itemize}
    then $\Aut(G)$ has a fixed point in $X$.
\end{thm}

Using the notion of the fixed-point dimension $\Fixdim$ (see Section \ref{sect 7}), we have:

\begin{cor}\label{main thm 2}
    Let $G=\ast_{i=1}^m \ast_{j=1}^{\ell_i} G_{i,j}$ be a free product as in Eq. (\ref{free decomposition}) satisfying the following properties:
    \begin{itemize}
        \item each $\ell_i\geq 2d+2$;
        \item $\min\{\Fixdim(G_{i,j}), ~\Fixdim(\Aut(G_{i,j}))\}\geq d$
    for each $G_{i,j}$. 
    \end{itemize}
    Then $\Fixdim(\Aut(G))\geq d$.
\end{cor}

Since the fixed-point dimensions of any finite group and its automorphism group are both infinite, we have a direct application of Corollary \ref{main thm 2} as follows.

\begin{cor}\label{cor of finite groups}
    Let $G=\ast_{i=1}^m \ast_{j=1}^{\ell_i} G_{i,j}$ be a free product of finite groups such that $G_{i_1,j_1}\cong G_{i_2,j_2}$ if and only if $i_1=i_2$. If each $\ell_i\geq n$, then $\Fixdim(\Aut(G))\geq \lfloor\frac{n}{2}\rfloor-1$.
\end{cor}

Note that Varghese \cite[Theorem A]{OV(Coxeter)} showed: $\Fixdim(\Aut(W_n)) \geq \lfloor\frac{n}{2}\rfloor-1$ for the universal Coxeter group $W_n=\ast_{n}\Z/2$. Corollary \ref{cor of finite groups} generalizes this result to the free product of general finite groups.

Furthermore, to prevent Properties (I) and (II) in Theorem \ref{main thm} from relying on the action on the space, we introduce properties determined only by $G$. Section \ref{sect 3} characterizes the relation between Property (II) and the first cohomology of $\Aut(G_{i,j})$ while Section \ref{sect 6} replaces Property (I) with a stronger assumption, Property (III). 

\begin{thm}\label{main2}
    Let $G=\ast_{i=1}^m \ast _{j=1}^{\ell_i}G_{i,j}$ be a free product as in Eq. (\ref{free decomposition}) and let $\Aut(G)$ act by isometries on a complete $\CAT0$ space $X$ with $\dim(X)\leq d$. If the following hold:
    \begin{itemize}
        \item each $\ell_i\geq 2d+2$;
        \item each $G_{i,j}$ satisfies Property (III);
        \item $\Aut(G)$ is finitely generated and $H^{1}(\Aut(G),\Z)=0$,
    \end{itemize}
    then $\Aut(G)$ has a fixed point in $X$.
\end{thm}
Note that the results above do not consider free products with free parts. A stronger conclusion (\textbf{Theorem \ref{main thm 3}}) also holds for the case with free parts, which recovers Theorem \ref{Bridson theorem} and Theorem \ref{main2}. However, since the process is rather cumbersome, we place it in Appendix \ref{section general}. To prove Theorem \ref{main thm} and Theorem \ref{main thm 3}, we introduce a technical tool, the $\FR$-graph (see Section \ref{sect 4} and Appendix \ref{section general}). Its construction arises from the Fouxe–Rabinovitch subgroups of $\Aut(G)$.

Finally, applying Theorem \ref{main2} to the free products of closed orientable surface groups, we have:

\begin{cor}\label{surface fixed point}
    Let $$G=\overbrace{\pi_1(\Sigma_{n_1})\ast\cdots\ast\pi_1(\Sigma_{n_1})}^{\ell_1}\ast
    \cdots \ast \overbrace{\pi_1(\Sigma_{n_m})\ast\cdots\ast \pi_1(\Sigma_{n_m})}^{\ell_m}\quad (n_i\neq n_j\geq 1 ~{for}~ i\neq j).$$
    Then any isometric action of $\Aut(G)$ on a complete $\CAT0$ space $X$ with $\dim(X)\leq d$ has a fixed point if each $\ell_i\geq 2d+2$.
\end{cor}

The paper is organized as follows. In Section \ref{sect 2}, we give some basic notations, definitions and introduce Bridson's bootstrap technique. In Section \ref{sect 3}, we introduce two properties (Properties (I) and (II)) and prove a key lemma which is the basis of Theorem \ref{main thm}. In Section \ref{sect 4}, we introduce two important tools, the $\FR$-graph and the rooted tree. In Section \ref{sect 5}, we finish the proof of Theorem \ref{main thm}. In Section \ref{sect 6}, we introduce $\phi$-torsion and Property (III). In Section \ref{sect 7}, we prove Theorem \ref{main2} and give some applications. Finally in Appendix \ref{section general}, we discuss the general cases and establish a result similar to Theorem \ref{main thm}, namely, Theorem \ref{main thm 3}.

\section{Preliminaries}\label{sect 2}

 In this section, we review some important facts on automorphisms of free products and group actions on $\CAT0$ spaces, see \cite{Andrew,Bridson(Hellytype),Bridson(NPC)} for further details.
 
\subsection{Automorphisms of free products}\label{sect 2.1}
For a finitely generated group $G$, the following decomposition is well known.

\begin{thm}[Grushko decomposition]
    Any finitely generated group $G$ can be decomposed as a free product:
    $G=G_1 \ast G_2 \ast \dots \ast G_\ell \ast F_r,$
    where each free factor $G_i\ncong \Z$, is nontrivial, freely indecomposable and $F_r$ is a free group of rank $r$.
\end{thm}

We now suppose that $G$ is such a free product with no free part $F_r$,
and suppose that $G$ is finitely generated with generators 
$$\s_G=\sqcup_{i=1}^{\ell} \s_{G_i},$$ 
where $\s_{G_i}$ is a finite generating set of $G_{i}$ such that $\phi(\s_{G_i})=\s_{\phi(G_i)}$ for any permutation $\phi$ between free factors in the same isomorphism class. After classifying free factors by (isomorphism) class, we rewrite this decomposition by
\begin{align}\label{free decomposition}
G=G_{1,1}\ast\cdots\ast G_{1,\ell_1}\ast \cdots \ast G_{m,1}\ast\cdots\ast G_{m,\ell_m}=\ast_{i=1}^m\ast_{j=1}^{\ell_i} G_{i,j},
\end{align}
where $G_{i_1,j_1}\cong G_{i_2,j_2}$ if and only if $i_1=i_2$. Following the notation in \cite{Andrew}, we consider three distinguished subgroups of $\Aut(G)$:
\begin{itemize}
    \item $\Fact(G)$: generated by automorphisms of single free factors.
    \item $\Perm(G)$: generated by permutations of free factors in the same class.
    \item $\FR(G)$: Fouxe–Rabinovitch subgroup.
\end{itemize} 
The above $\FR(G)$ can be finitely generated by \emph{partial conjugations} $\PC(G)$, where
    \begin{align}\label{PC}
        \PC(G): =\{(A,b)\in \Aut(G)\mid \text{$A, B$ are distinct free factors of $G$ in Eq. (\ref{free decomposition}) and $b\in \s_B$}\},
    \end{align}
and 
$$(A,b)(g): =\begin{cases}
    b^{-1}gb\quad & g\in \s_A\\
    g\quad & g\in \s_G-\s_A
\end{cases}.$$  
We will denote $(A,B): =\{(A,b)\mid b\in B\}\subset \FR(G)$ for convenience.

Note that $(A, b)(A, b')=(A, bb')$ for any $b, b'\in B$. Moreover, we have:

\begin{lem}\label{conjugate}
Let $A,B,C,D$ be four distinct free factors of $G$. The following hold.
\begin{enumerate}
    \item For any $\phi\in \Perm(G)$ and $b\in B$, we have $\phi(A,b)\phi^{-1}=(\phi(A),\phi(b))$ and hence $$\phi(A,B)\phi^{-1}=(\phi(A),\phi(B))\subset \FR(G).$$
    \item There is an embedding $$B\hookrightarrow \Aut(G), \qquad b \mapsto (A,b),$$ which induces an isomorphism $(A, B)\cong B$. Hence $(A, B)$ is finitely generated. 
    \item For any $b\in B$ and $d\in B\cup D$, we have 
    $$(C,d)^{-1}(A,b)(C,d)=(A,b),$$ and hence the subgroups $(A,B)$ and $(C,B)$ commute, so do $(A,B)$ and $(C,D)$.
    \item For any $b\in B$ and $c\in C$, we have $$(B,c)^{-1}(A,b)(B,c)=(A,c)(A,b)(A,c)^{-1},$$ and hence $(B,C)$ normalizes $\generate{(A,B)\cup (A,C)}$.
\end{enumerate}
\end{lem}
\begin{proof}
(1) For any $b \in B$ and $g\in \s_G$,
    $$\phi(A,b)\phi^{-1}(g)=\begin{cases}
        \phi(b^{-1}\phi^{-1}(g)b)=\phi(b)^{-1}g\phi(b) \quad & \phi^{-1}(g)\in \s_{A}\\
        \phi(\phi^{-1}(g))=g \quad & \phi^{-1}(g)\in \s_{G}-\s_{A}
    \end{cases}.$$

(2), (3) and (4) can be verified by definition directly.
\end{proof}

We reformulate \cite[Proposition 3.1]{Gilbert(presentation)} (cf. \cite[Proposition 2.8]{Andrew}) and \cite[Proposition 2.9]{Andrew} in the following form required for our arguments:

\begin{thm}\label{Gilbert}
    Let $G=\ast_{i=1}^m\ast_{j=1}^{\ell_i} G_{i,j}$ be a free product as in Eq. (\ref{free decomposition}) (note that $G_{i,j}\ncong F_r$). Then
    \begin{itemize}
        \item $\Aut(G)=\FR(G)\rtimes \generate{\Fact(G), \Perm(G)}$;
        \item $\generate{\Fact(G), \Perm(G)}=\Pi_{i=1}^{m}(\Aut(G_{i,1})^{\ell_i}\rtimes S_{\ell_i})$.
    \end{itemize}
\end{thm}

\subsection{\texorpdfstring{$\rm CAT(0)$ spaces}{CAT(0) spaces}} 
Let $X$ be a geodesic metric space and let $\Isom(X)$ denote the isometry group of $X$.
For any geodesic triangle $\Delta\subset X$, there is a triangle $\overline{\Delta}\subset \mathbb E^2$ with the same edge lengths as $\Delta$. Let $\bar{x}\mapsto x$ denote the map $\overline{\Delta}\to \Delta$ that sends each side of $\overline \Delta$ isometrically onto the corresponding side of $\Delta$. One says that $X$ is a \emph{$\CAT0$ space} if for all $\Delta$ and all $\bar x,\bar y\in \overline{\Delta}$, the inequality 
$$d_X(x,y)\leq d_{\mathbb E^2}(\bar x,\bar y)$$
holds. A \emph{complete $\CAT0$ space} means that the space is complete under the geodesic metric. 
In a complete $\CAT0$ space $X$, there is a unique geodesic $[x, y]$ joining each pair of points $x, y\in X$. A subspace $Y\subset X$ is said to be \emph{convex} if $[y,y']\subset Y$ whenever $y,y'\in Y$.

Given a subset $H\subset \Isom(X)$, we denote its fixed point set
$$\Fix(H):=\{x \in X \mid  h(x)=x,~\forall h \in H\}.$$ 
Note that $\Fix(H)=\Fix(\generate{H})$ is closed and convex, and hence a complete $\CAT0$ subspace of $X$. Moreover, given two subsets $H_1,H_2 \subset \Isom(X)$ and an element $g\in \Isom(X)$, there are two fundamental facts:
\begin{itemize}
    \item $\Fix(\{H_1,H_2\})=\Fix(H_1\cup H_2)=\Fix (H_1)\cap \Fix (H_2)$;
    \item $\Fix(H_2)=g\Fix(H_1)$, if $\generate{H_2}=g\generate{H_1}g^{-1}$.
\end{itemize}

A classical result of Bruhat and Tits states that any action of a finite group on a complete $\CAT0$ space has a fixed point (see \cite[Corollary II.2.8(1)]{Bridson(NPC)}). Moreover, Bridson obtained:
\begin{prop}\cite[Proposition 1.1]{Bridson(Hellytype)}\label{Bruhat and Tits}
    Given a complete $\CAT0$ space $X$ and any isometric group action $\Gamma$, the following are equivalent:
\begin{itemize}
    \item $\Gamma$ has a bounded orbit in $X$;
    \item all $\Gamma$-orbits are bounded;
    \item $\Fix(\Gamma)$ is nonempty.
\end{itemize}
\end{prop}

\subsection{Nerves of subsets}
Let $X$ be a topological space and $\C$ be a collection of subsets of $X$. The \textit{nerve} of $C$, denoted by $\N(\C)$, is the simplicial complex whose $k$-simplices $[i_0, \ldots, i_k]$ correspond to sub-collections $\{C_{i_0}, \ldots, C_{i_k}\}\subset C$ with nonempty intersection.

\begin{thm}\cite[Theorem B]{Bridson(Hellytype)}\label{Nerve}
    Let $X$ be a complete convex metric space and let $\C$ be a finite collection of closed convex subsets of $X$. If $\dim(X)\leq d$, then every continuous map $\N(\C)\to \mathbb S^r$ from the nerve of $\C$ to a sphere of dimension $r\geq d$ is homotopic to a constant map.
\end{thm}

A useful corollary of Theorem \ref{Nerve} is the following metric Helly-type result.

\begin{cor}\cite[Corollary 2.4]{Bridson(Hellytype)} \label{Helly type theorem}
Let $C_1,\ldots,C_m$ be closed convex subspaces in a convex metric space $X$. If $\dim(X)\leq d$ and $\cap_{i\in I} C_{i}\neq \emptyset$ for each $I\subset \{1,\ldots,m\}$ with $|I|\leq d+1$, then $\cap_{i=1}^m C_i \neq \emptyset$.
\end{cor}

\subsection{Normalized and commuting subgroups, bootstrap and conjugate bootstrap}
\begin{prop}\cite[Proposition 1.3]{Bridson(Hellytype)}\label{Normalize maintain fixed point}
    Let $G$ be a group, let $H_1,H_2<G$ be subgroups, and suppose that $H_2$ normalizes $H_1$. Whenever $G$ acts by isometries on a complete $\CAT0$ space, if $\Fix(H_1)$ and $\Fix(H_2)$ are nonempty, then $\Fix(H_1)\cap \Fix(H_2)$ is nonempty.
\end{prop}

As a direct corollary, we have:

\begin{cor}\label{product maintain fix}
    Let $X$ be a complete $\CAT0$ space and let $H_1,\ldots,H_\ell$ be subgroups of $\Isom(X)$ with $\Fix(H_i)$ nonempty for $i = 1,\ldots,\ell$. Then 
    \begin{itemize}
        \item whenever $H_1,\ldots,H_\ell$ commute pairwise, $\bigcap_{i=1}^{\ell}\Fix(H_i)$ is nonempty;
        \item whenever $\generate{H_1,H_2}=H_1\rtimes H_2$, $\Fix(H_1)\cap \Fix(H_2)$ is nonempty.
    \end{itemize}
\end{cor}

\begin{lem}\cite[Lemma 1.5]{Bridson(Hellytype)}\label{quotient action induced by group}
    Let $N<G$ be a subgroup of finite index. If $G$ acts by isometries on a complete $\CAT0$ space and $\Fix(N)$ is nonempty, then $\Fix(G)$ is nonempty.
\end{lem}

The following bootstrap technique can be found in Bridson \cite[Proposition C, Corollary 4.6]{Bridson(Hellytype)}.

\begin{prop}[\cite{Bridson(Hellytype)}, Bootstrap Lemma]\label{Bootstrap}
    Let $k_1,\ldots,k_n$ be positive integers and let $X$ be a complete $\CAT0$ space with $\dim(X)<k_1+\cdots +k_n$. Let $S_1,\ldots,S_n\subset \Isom(X)$ be subsets with $[s_i,s_j]=1$ for all $s_i\in S_i$ and $s_j\in S_j(i\neq j)$.
    
    If, for $i=1,\ldots,n$, each $k_i$-element subset of $S_i$ has a fixed point in $X$, then for some $i$, every finite subset of $S_i$ has a fixed point.
\end{prop}

\begin{cor}[\cite{Bridson(Hellytype)}, Conjugate Bootstrap]\label{Conjugate Bootstrap}
Let $k$ and $n$ be positive integers and let $X$ be a complete $\CAT0$ space with $\dim(X)<nk$. Let $S_1,\ldots,S_n$ be commuting (i.e. $[s_i,s_j]=1$ for all $s_i \in S_i$ and $s_j \in S_j (i \neq j)$) conjugates of a subset $S \subset \Isom(X)$.

If each $k$-element subset of $S$ has a fixed point in $X$, then so does each finite subset of $S$ and of $S_1 \cup \ldots \cup S_n$. Additionally, if $S<\Isom(X)$ is a finitely generated subgroup, then $\Fix(S)\neq \emptyset$. 
\end{cor}

Given a group $G$ and the symmetric group $S_{d}$ of degree $d$, we define the wreath product throughout this paper by
$$G \wr S_{d} := \bigoplus_{1\leq i\leq d} G_i \rtimes S_{d},$$ where each $G_i$ is a copy of $G$ and $S_d$ acts by permuting the indices. There is a natural embedding $G \hookrightarrow G \wr S_{d}$ through $G_1$.
By using results above, we have the following key lemma.
\begin{lem}\label{close under conjugation, radical, multiplication}
    Let $G\wr S_{d+1}$ act by isometries on a complete $\CAT0$ space $X$ with $\dim(X)\leq d$.
    Then
    $$\Elp_X(G):=\{g\in G\mid \text{$\Fix (g)\neq \emptyset$}\}$$ is closed under conjugation, taking roots, and multiplication. Moreover, every finite subset of $\Elp_X(G)$ has a fixed point.
\end{lem}
\begin{proof}
    $\Elp_X(G)$ is closed under conjugation clearly. Take any $g_1,g_2\in \Elp_X(G)$ and $x\in \Fix(g_1)$. If $h$ is a root of $g_1$, i.e., $h^m=g_1$, then the orbit $\generate{h}x$ is finite and thus $\Fix(h)\neq \emptyset$ by Proposition \ref{Bruhat and Tits}. Furthermore, by hypothesis, we can choose $d+1$ copies of $G$ which are pairwise commuting and conjugate. Thus, by the conjugate bootstrap (Corollary \ref{Conjugate Bootstrap}), $\{g_1,g_2\}$ has a fixed point and so does $g_1g_2\in \generate{g_1,g_2}$. The same choice of commuting conjugates implies that any finite subset of $\Elp_X(G)$ has a fixed point.
\end{proof}

\section{Properties (I)(II) and a key lemma}\label{sect 3}

Let $G=\ast_{i=1}^m\ast_{j=1}^{\ell_i} G_{i,j}$ be a free product as in Eq. (\ref{free decomposition}), and let $\Aut(G)$ act by isometries on a complete $\CAT0$ space $X$. Note that for every two distinct free factors $A$ and $B$, by Lemma \ref{conjugate}, there is an embedding $$B\hookrightarrow \Aut(G), \qquad b \mapsto (A,b).$$
Therefore, we have an isomorphism $(A, B)\cong B$. 
\subsection{Properties (I)(II)}
Under the action of $\Aut(G)$ on $X$, we introduce two properties:
\begin{itemize}
    \item \emph{Property (I)}: every subgroup $(A, B)$ has a fixed point in $X$.
    \item \emph{Property (II)}: every subgroup $\Aut(A)$ has a fixed point in $X$.
\end{itemize}

Since $\Perm(G)$ can permute the isomorphic but distinct free factors of $G$, Property (I) holds if and only if  
both $(G_{i,1},G_{i,2})$ and $(G_{i,1},G_{i',1})$ have a fixed point in $X$ for each $1\leq i\neq i'\leq m$; Property (II) holds if and only if $\Aut(G_{j,1})$ has a fixed point for each $1\leq j\leq m$.

To give more formulation of Property (II) in Proposition \ref{property 2 equal finite abelianization}, we first prove the following proposition.
\begin{prop}\label{wreath product}
    Let $S_\ell$ be the symmetric group of degree $\ell$ and $H$ be a finitely generated group. Then the first cohomology $H^1(H,\Z)=0$ if and only if the wreath product $$W_\ell=H\wr S_\ell$$
    has a fixed point (and so does $H$) whenever $W_\ell$ acts by isometries on a complete $\CAT0$ space $X$ with $\dim(X)=d\leq \lfloor\frac{\ell}{2}\rfloor-1$.
\end{prop}
\begin{proof}
    The ``only if'' part.
    Consider $$W_2=(H\times H)\rtimes \generate{\sigma},$$
    where $\sigma$ is the involution interchanging the two factors and $K=\langle\!\langle\sigma\rangle\!\rangle_{W_2}$ is the normal closure of $\sigma$ in $W_2$.
    We claim that $$ W_2/K\cong H_{\mathrm{ab}}$$ is finite, where $H_{ab}$ is the abelianization of $H$.
    Indeed, the epimorphism
    \begin{equation}\label{beta}
        \beta:W_2\twoheadrightarrow H_{\mathrm{ab}}, \qquad \beta(ab\sigma^\varepsilon)=[a]+[b]
    \end{equation}
    is well-defined, and its kernel is $K$, where $[a], [b]$ denote the images of $a, b$ in $H_{ab}$. Moreover, $H^1(H,\Z)=0$ implies $H_{\mathrm{ab}}$ is finite. 

    Since $W_2$ is finitely generated and $W_2/K$ is finite, $K$ is finitely generated.
    By hypothesis, $\ell\geq 2d+2$, we can choose $d+1$ pairwise commuting and conjugate copies of $W_2$ in $W_\ell$ and have an embedding
    $$W_2\wr S_{d+1}\hookrightarrow W_{\ell}$$
    defined by
    \begin{equation*}
        \begin{aligned}
            S_{d+1} \hookrightarrow S_\ell&: \quad  (i,j) \mapsto (2i-1,2j-1)(2i,2j),\\
            \bigoplus_{d+1}W_2\hookrightarrow W_\ell&: \quad a_1b_1\sigma^{\epsilon_1} \cdots a_{d+1}b_{d+1}\sigma^{\epsilon_{d+1}} \mapsto a_1b_1\cdots  a_{d+1}b_{d+1}(1,2)^{\epsilon_{1}}\cdots(2d+1,2d+2)^{\epsilon_{d+1}}.
        \end{aligned}
    \end{equation*}
    Since $\sigma\in \Elp_X(W_2)$, by Lemma \ref{close under conjugation, radical, multiplication}, $K\subset\Elp_X(W_2)$, and thus the finite generating set of $K$ has a fixed point. Then $W_2$ has a fixed point by the above claim and Lemma \ref{quotient action induced by group}, and so does $H$. Therefore $W_{\ell}$ has a fixed point by Corollary \ref{product maintain fix}.

    The ``if'' part. Take a similar epimorphism as in Eq. (\ref{beta}): $$\beta: W_\ell \twoheadrightarrow H_{\rm ab}.$$
    If $H^1(H,\Z)\neq 0$, then we have epimorphisms $$\psi:H_{\rm ab}\twoheadrightarrow \Z,\quad \psi \circ \beta:W_{\ell}\twoheadrightarrow \Z.$$
    However, $\Z$ has a translation action on the line with no fixed point and so does $W_\ell$, a contradiction.
\end{proof}

\begin{remark}
    Proposition \ref{wreath product} is a similar higher-dimensional version of \cite[Theorem 1.1]{CK}. Our result extends this fixed-point phenomenon from actions on trees to isometric actions on complete $\CAT0$ spaces of finite dimension.
\end{remark}

\begin{prop}\label{property 2 equal finite abelianization}
    Let $G=\ast_{i=1}^m \ast _{j=1}^{\ell_i}G_{i,j}$ be a free product as in Eq. (\ref{free decomposition}) with $\Aut(G)$ finitely generated, and let $\Aut(G)$ act by isometries on a complete $\CAT0$ space $X$ with $\dim(X)\leq d$.
    \begin{enumerate}
        \item If $\ell_i \geq d+1$, then $\Aut(G_{i,j})$ has a fixed point if and only if $\Aut(G_{i,j})$ has a finite generating set $\s_{\Aut(G_{i,j})}$ such that every generator in $\s_{\Aut(G_{i,j})}$ has a fixed point.
        \item If $\ell_i \geq 2d+2$ and $H^1(\Aut(G_{i,1}),\Z)=0$, then $\Aut(G_{i,j})$ has a fixed point.
    \end{enumerate}
\end{prop}
\begin{proof}
    Theorem \ref{Gilbert} showed that
    \begin{equation}\label{Gilbert form}
        \Aut(G)=\FR(G)\rtimes \Pi_{i=1}^{m}(\Aut(G_{i,1})^{\ell_i}\rtimes S_{\ell_i}).
    \end{equation}
    Thus $\Aut(G)$ is finitely generated if and only if every $\Aut(G_{i,j})$ is finitely generated.
    
    (1) It suffices to prove the ``if'' part. By hypothesis, every $\ell_i\geq d+1$, we can choose $d+1$ pairwise commuting and conjugate copies of $\Aut(G_{i,1})$, so there is an embedding 
    \begin{equation*}
        \Aut(G_{i,1})\wr S_{d+1}\hookrightarrow \Aut(G).
    \end{equation*}
    Since the finite generating set of $\Aut(G_{i,1})$ is a subset of $\Elp_X(\Aut(G_{i,1}))$, it has a fixed point by Lemma \ref{close under conjugation, radical, multiplication}.
       
    (2) By Eq. (\ref{Gilbert form}), $\Aut(G_{i,1})\wr S_{\ell_i}\hookrightarrow \Aut(G)$ acts by isometries on $X$. Then by Proposition \ref{wreath product}, $\Aut(G_{i,j})$ has a fixed point.
\end{proof}
When some $H^1(\Aut(G_{i,j}),\Z)\neq 0$, we have:
\begin{prop}
     Let $G=\ast_{i=1}^m \ast _{j=1}^{\ell_i}G_{i,j}$ be a free product as in Eq. (\ref{free decomposition}) such that $$H^1(\Aut(G_{i,j}),\Z)\neq 0$$ for some free factor $G_{i,j}$. Then $\Aut(G)$ can act by isometries on the line with no fixed point. 
\end{prop}
\begin{proof}
    By Eq. (\ref{Gilbert form}),
    we have epimorphisms
    $$\Aut(G)\twoheadrightarrow \Aut(G_{i,j})_{\rm ab} \twoheadrightarrow \Z$$ for some $H^1(\Aut(G_{i,j}),\Z)\neq 0$.
    Thus, $\Aut(G)$ has a translation action on the line with no fixed point.
\end{proof}

\subsection{A key lemma}

\begin{lem}\label{any two have common fixed point,high level}
    Let $G=\ast_{i=1}^{m} \ast_{j=1}^{\ell_i}G_{i,j}$ be a free product as in Eq. (\ref{free decomposition}) and let $\Aut(G)$ act by isometries on a complete $\CAT0$ space $X$ with $\dim(X)\leq d$. If each $\ell_i\geq 2d+2$ and $G$ satisfies Property (I), then $$\Fix(G_1, G_2) \cap \Fix(G_3, G_4)\neq \emptyset$$ for every four free factors $G_1\neq G_2$ and $G_3\neq G_4$.  
\end{lem}

\begin{proof}
     We consider all possible pairs $(G_1, G_2)$ and $(G_3, G_4)$. Let $A, B, C, D$ denote distinct free factors (possibly isomorphic). After removing repeating cases, there are exactly five cases as follows.
     
(1) $(A,B)$ and $(C,D)$. By Lemma \ref{conjugate}, these two subgroups commute, and hence by Corollary \ref{product maintain fix}, $$\Fix(A,B) \cap \Fix(C,D) \neq \emptyset.$$
      
      (2) $(A,B)$ and $(C,B)$. By the same arguments as in item (1), $\Fix(A,B) \cap \Fix(C,B) \neq \emptyset.$

    (3) $(A,B)$ and $(B,A)$. Since $\ell_i\geq 2d+2$, we can take the following $2d+2$ distinct free factors:
    $$A=A_0, A_1, \ldots, A_d, \quad B=B_0, B_1, \ldots, B_d,$$ where all the $A_k\cong A$ and $B_k\cong B$. Then we have $d+1$ subsets of $\FR(G)$:
   \begin{equation}\label{eq. S_k}
  S_k:=(A_{k}, B_{k})\cup (B_{k},A_{k}), \quad k=0,1,\ldots,d, 
   \end{equation}
     which are commuting and conjugate pairwise by Lemma \ref{conjugate}. Moreover, Property (I) implies that every element in $S_k$ has a fixed point. Then by the conjugate bootstrap (Corollary \ref{Conjugate Bootstrap}), each finite subset of $S_0$ has a fixed point. It follows that 
     $$\Fix (S_0)=\Fix(A, B) \cap \Fix(B,A)\neq \emptyset.$$

    (4) $(A,B)$ and $(A,D)$. In this case, we can take $$S_0=(A, B) \cup (A, D) \ \ \text{and}\ \  S_k=(A_{k}, B) \cup (A_{k}, D)$$ for $1\leq k \leq d$. Then following similar arguments as in item (3), we have $\Fix(A,B) \cap \Fix(A,D) \neq \emptyset$.

    (5) $(A,B)$ and $(B,D)$. 
    
    Whenever $\ell_i\geq 2d+3$, or $A, B, D$ are not in the same isomorphic class, by the same arguments as in item (3), but with Eq. (\ref{eq. S_k}) replaced by the following $d+1$ subsets of $\FR(G)$:
    $$S_k:=(A_{k}, B_{k})\cup (B_{k}, D),  \quad k=0,1,\ldots,d, $$ we conclude that $\Fix(A, B) \cap \Fix(B,D)\neq \emptyset.$

    Whenever $\ell_i=2d+2$ and $A\cong B\cong D$, we again claim that $\Fix(A,B)\cap\Fix(B,D)\neq \emptyset$. Otherwise, if $\Fix(A,B)\cap\Fix(B,D)=\emptyset$, then we take 
    $$S:=(A, B) \cup (B, D), \quad S':=(D, B) \cup (B, A),$$ 
    \begin{equation}\label{eq. CC}
      \C:=\{\Fix(A,B), ~\Fix(B, D), ~\Fix(D, B), ~\Fix(B, A)\}. 
    \end{equation}
    Since $S'$ is conjugate to $S$, we have $\Fix (S') =\Fix (S) = \emptyset$. 
    By the above items (2)--(4), we have
    $$\Fix(A,B) \cap \Fix(D,B) \neq \emptyset, \quad \Fix(B,A) \cap \Fix(B,D) \neq \emptyset;$$
    $$\Fix(A,B) \cap \Fix(B,A) \neq \emptyset, \quad \Fix(D,B) \cap \Fix(B,D) \neq \emptyset.$$
    Therefore, the nerve $\N(\C)$ is a square, see Figure \ref{cube}. 
     \tikzset{every picture/.style={line width=0.75pt}} %set default line width to 0.75pt
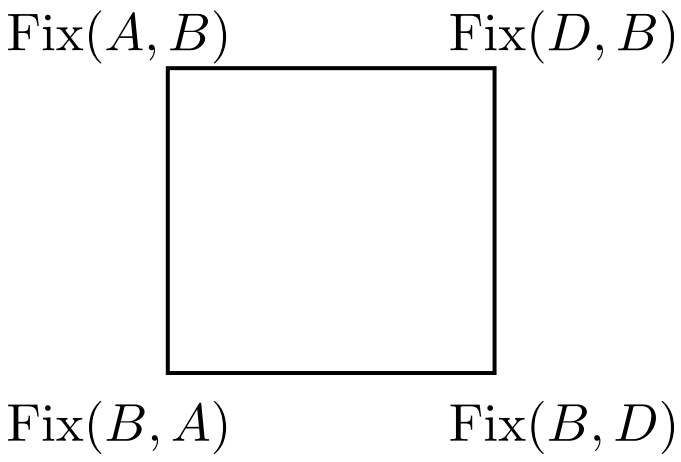
\begin{figure}[htbp]
\begin{tikzpicture}[x=0.75pt,y=0.75pt,yscale=-1,xscale=1]

\draw   (241,111) -- (323,111) -- (323,187.5) -- (241,187.5) -- cycle ;

% Text Node
\draw (199,95) node [anchor=north west][inner sep=0.75pt]   [align=left] {$\Fix(A,B)$};
% Text Node
\draw (199,193) node [anchor=north west][inner sep=0.75pt]   [align=left] {$\Fix(B,A)$};
% Text Node
\draw (310,95) node [anchor=north west][inner sep=0.75pt]   [align=left] {$\Fix(D,B)$};
% Text Node
\draw (310,193) node [anchor=north west][inner sep=0.75pt]   [align=left] {$\Fix(B,D)$};
\end{tikzpicture}
\caption{Nerve of $\C$ in Eq. (\ref{eq. CC})}\label{cube}
\end{figure}

Moreover, there are $3\cdot  \lfloor \frac{2d+2}{3} \rfloor$ distinct free factors:  $A_{k}, B_{k}, D_k$, $k=1, \ldots,  \lfloor \frac{2d+2}{3} \rfloor$. We now take 
    $$\Tilde S_k:=(A_{k}, B_{k})\cup (B_{k}, D_{k})\cup(B_{k},A_{k})\cup (D_{k},B_{k})$$ and
    $$\C_k:=\{\Fix(A_{k}, B_{k}), ~\Fix(B_{k}, D_{k}), ~\Fix(B_{k},A_{k}), ~\Fix(D_{k},B_{k})\}.$$ 
    Then all of $\tilde S_k$ are conjugate to $S\cup S'$ and hence $\N(\C_k)\cong \N(\C)$ are squares. Since elements in different $\tilde S_k$ commute pairwise, by \cite[Lemma 4.1]{Bridson(Hellytype)},
    \begin{equation*}
             \N(\cup_{i}\C_i)=\N(\C_1)\ast\N(\C_2)\ast \cdots \ast \N(\C_{\lfloor \frac{2d+2}{3} \rfloor})=\mathbb S^{2\cdot\lfloor \frac{2d+2}{3} \rfloor-1}.
    \end{equation*}
    That is, $\N(\cup_{i}\C_i)$ is a sphere of dimension $2\cdot\lfloor \frac{2d+2}{3} \rfloor-1 \geq  d$, which contradicts Theorem \ref{Nerve}.
\end{proof}

\section{\texorpdfstring{The $\FR$-graph $\G(S)$ and rooted tree $T_v^\G$}{FR-graph}}\label{sect 4}

Recall that for $G=\ast_{i=1}^m\ast_{j=1}^{\ell_i} G_{i,j}$ in Eq. (\ref{free decomposition}), the Fouxe–Rabinovitch subgroup $\FR(G)\subset \Aut(G)$ is finitely generated by partial conjugations. Let 
\begin{equation}\label{eq. FR0}
  \FR_0(G):=\{(A,b)\in (A,B)\mid \text{$A, B$ are distinct free factors of $G$}\}\subset\FR(G).  
\end{equation}
Note that $\langle \FR_0(G) \rangle=\FR(G)$.
\subsection{\texorpdfstring{$\FR$-graph $\G(S)$}{FR-graphs}}
For any subset $S\subset \FR_0(G)$ ({\bf since any action of $\id$ is trivial, we replace $S$ with $S\setminus\id$ if $\id \in S$} in the following), 
we define three finite sets as follows:
\begin{align*}
  \Cf(S)&:=\{A\mid \text{$(A,b)\in S$}\},\\
\V(S)&: =\{A \mid (A,b)\in S \text{, or } (B,a) \in S \text{ for some }a \in A \},\\
\E(S)&: = \{(A,B)\mid (A,b)\in S \text{ for some $(A,b)\in (A,B)$}\}.
\end{align*} 
Note that the set $\V(S)$ consists of two types of free factors: all free factors occurring in the first entry of each pair in $S$ (i.e., conjugated free factors $\CF$), together with those free factors containing elements from the second entry of such pairs. The set $\E(S)$ denotes a family of subgroups of $\FR(G)$ that contain some elements of $S$. Clearly, we have
\begin{equation}\label{eq.FVE}
|\E(S)|\leq |S|\quad \mbox{and} \quad |\CF(S)|\leq \min\{|\V(S)|, |\E(S)|\}.
\end{equation}
If $(A,B)\in \E(S)$, then $A,B\in \V(S)$. Thus we can construct a directed graph $\G(S)$ as follows:
\begin{defn}[$\FR$-graph]
    For any subset $S\subset \FR_0(G)$, the \emph{$\FR$-graph generated by $S$}, denoted by $\G(S)$, is a finite directed graph with
    \begin{itemize}
        \item vertex set $\Vcal(\G(S))$ : for each $A \in \V(S)$, there is a vertex $v_A$ with the \emph{vertex group} $G_{v_A}:=A$;
        \item edge set $\Ecal(\G(S))$: for each $(A,B) \in \E(S)$, there is a directed edge $e_{A,B}$ from $v_A$ to $v_B$ with the \emph{edge group} $G_{e_{A,B}}:=(A,B)$.
    \end{itemize}
    
\end{defn}

Note that an $\FR$-graph is a directed graph whose vertices and edges are each assigned a corresponding group. Specifically, its vertex groups are free factors of $G$, while its edge groups are subgroups of $\FR(G)$ of the form $(A,B)$, where $A$ and $B$ correspond to the origin and terminus of the edge, respectively.
Since $(A,B)$ and $(B,A)$ are two distinct subgroups of $\FR(G)$, there exist at most two directed edges between any pair of distinct vertices, one in each direction, and no loop can exist.
For any $\FR$-graph $\G$, we denote the union of all its edge groups $G_e$ as
$$\M(\G):=\cup_{e\in \Ecal(\G)} G_e.$$ It is immediate that $\G=\G(\M(\G))$. Furthermore, we have the following straightforward observations.

\begin{lem}\label{basic lemma of FR-graph} 
Let $S, S'$ be two subsets of $\FR_0(G)$. Then the following hold.
\begin{enumerate} 
    \item  $S\subset \M(\G(S))$, ~~$\CF(S)=\CF(\M(\G(S)))\subset\V(S)=\V(\M(\G(S)))$ and $\E(S)=\E(\M(\G(S)))$.
    
    \item If $S \neq \emptyset$, then each connected component of $\G(S)$ possesses at least two vertices and one edge.
    
    \item  If $S\subset S'$, then $\G(S)$ is a subgraph of $\G(S')$, with all vertex and edge groups inherited. Moreover, $\M(\G)$ is the union of subsets $S$ with $\G(S)=\G$.
     
    \item  Every $\FR$-graph is a subgraph of $\G(\FR_0)$, which possesses $\sum_{i=1}^m \ell_i$ vertices and admits two directed edges between every pair of distinct vertices.
\end{enumerate}
\end{lem}
Note that there may exist infinitely many distinct sets $S$, and distinct choices of $S$ can produce the same graph $\G(S)$. However, by Lemma \ref{basic lemma of FR-graph}(4), there are only finitely many distinct graphs $\G(S)$. For this reason, it is much more convenient to study $\G(S)$ instead of $S$ itself directly.

\subsection{\texorpdfstring{Rooted tree $T_v^\G$}{Rooted trees}}
Let $S\subset \FR_0(G)$ and $\G:=\G(S)$. For an edge $e$ in $\G$, let $o(e)$ and $\tau(e)$ denote the origin and terminus of $e$, respectively. 
A \emph{directed path} from vertex $v$ to vertex $u$ is a sequence of directed edges $e_1,\ldots, e_n (n\geq 1)$ such that $o(e_1)=v$, $\tau(e_n)=u$ and $\tau(e_i)=o(e_{i+1})$ for $i=1,\ldots,n-1$. Note that $u$ and $v$ may be equal, and in this case the directed path becomes a \emph{directed cycle}. Since there are no loops, every directed cycle consists of at least two directed edges. For a vertex $v$ in $\G$, we define the \emph{post-vertex set} $\Post^\G(v)$ as 
$$\Post^\G(v):=\{u\in \Vcal(\G)\mid \text{there exists a directed path from $v$ to $u$ in $\G$}\}.$$
\begin{defn}[Rooted tree $T^\G_v$]\label{Rooted trees}
For any vertex $v$ in $\G$ satisfying $v\notin \Post^\G(v)$, we construct a level-1 \emph{rooted tree} $T^\G_v$ as follows:
\begin{align*}
\text{ vertex set } \Vcal(T^\G_v) &:=\{v\}\cup \Post^\G(v),\\
\text{ edge set } \Ecal(T^\G_v) &:=\{(v,u) \mid u\in \Post^{\G}(v)\}.
\end{align*}
Note that $T_v^{\G}$ can be regarded as a special $\FR$-graph induced from $\G$: $\Vcal(T_v^{\G})$ is a subset of $\Vcal(\G)$. Hence, the vertex group $G_u$ of any vertex $u \in \Vcal(T_v^{\G})$ is exactly the vertex group of $u \in \Vcal(\G)$. For an edge $e=(v,u)$ in $\Ecal(T_v^{\G})$, the corresponding edge group $G_e$ is defined as $(G_v,G_u)$.
\end{defn}

\section{Proof of Theorem \ref{main thm}}\label{sect 5}
 
In this section, we always assume that $G=\ast_{i=1}^{m} \ast_{j=1}^{\ell_i}G_{i,j}$ is a free product as in Eq. (\ref{free decomposition}), and that $\Aut(G)$ acts by isometries on a complete $\CAT0$ space $X$. 

If $G$ satisfies Property (I), then every subgroup $(G_{i,j},G_{i',j'})$ of $\FR(G)$ has a fixed point in $X$. If, in addition, $\dim(X)=1$ and $\ell_i \geq 4$, then combining Lemma \ref{any two have common fixed point,high level} with Corollary \ref{Helly type theorem}, we obtain that $\FR(G)$ has a fixed point in $X$.

In general, we shall prove the following result.

\begin{prop}\label{solvable graph lemma} 
    Suppose $2\leq\dim(X)\leq d$, each $\ell_i\geq 2d+2$ and $G$ satisfies Property (I). Then for any subset $S\subset \FR_0(G)$, $\M(\G(S))$ has a fixed point in $X$. In particular, $\FR(G)=\generate{\M(\G(\FR_0))}$ has a fixed point in $X$.
\end{prop}

\begin{proof}[\textbf{Proof of Theorem \ref{main thm}}]

By Theorem \ref{Gilbert}, we have $\Aut(G)=\FR(G)\rtimes \generate{\Fact(G),\Perm(G)},$ and $$\generate{\Fact(G),\Perm(G)}=\Pi_{i=1}^{m}(\Aut(G_{i,1})^{\ell_i}\rtimes S_{\ell_i}).$$
    Since $G$ satisfies Property (II), i.e., each $\Aut(G_{i,1})$ has a fixed point, $\generate{\Fact(G),\Perm(G)}$ has a fixed point by Corollary \ref{product maintain fix}.
    Furthermore, by Proposition \ref{solvable graph lemma}, $\FR(G)$ has a fixed point in $X$. Hence $\Aut(G)$ has a fixed point in $X$ again by Corollary \ref{product maintain fix}.
\end{proof}

In the remainder of this section, we prove Proposition \ref{solvable graph lemma}. By disregarding the edge directions, i.e. treating it as an undirected graph, we first discuss whether $\G(S)$ is a tree or is connected.

If $|\Vcal(\G(S))|=|\Ecal(\G(S))|+1$ and $\G(S)$ is connected, then $\G(S)$ is a tree. If $|\Vcal(\G(S))|\geq |\Ecal(\G(S))|+2$, then $\G(S)$ must be disconnected. Accordingly, $\G(S)$ can be decomposed as $\G(S)=\G_1\sqcup \G_2$ where $\G_1$ and $\G_2$ are two nontrivial disjoint subgraphs (which are not necessarily connected components). Since subgroups $(A,B)$ and $(C,D)$ commute whenever $A,B,C,D$ are all distinct, the edge groups $G_e$ and $G_{e'}$ commute for all edges $e\in \G_1$ and $e'\in \G_2$.  Consequently, the set $S$ decomposes into $S=S_1\sqcup S_2$ such that $\G(S_i)=\G_i$. By the above argument, $S_1$, $S_2$ commute and so do $\M(\G_1)$ and $\M(\G_2)$.

Furthermore, we need several lemmas. First, by a simple calculation, we have:

\begin{lem}\label{ad function2}
    For any integers $d\geq k\geq 3$, we have $\lfloor\frac{2d+1}{k}\rfloor(k-1)>d.$
\end{lem}

We begin by considering the case $|\E(S)|\leq d$.

\begin{lem}\label{|E(S)|leq d}
    If $S$ in Proposition \ref{solvable graph lemma} satisfies $|\E(S)|\leq d$,  then $\M(\G(S))$ has a fixed point in $X$. 
\end{lem}
\begin{proof}
    We proceed by induction on $|\E(S)|$. When $|\E(S)|=1$, the conclusion is equivalent to Property (I); and when $|\E(S)|=2$, the existence of a fixed point for $\M(\G(S))$ follows from Lemma \ref{any two have common fixed point,high level}. We now assume that the conclusion holds for all subsets $S\subset \FR_0(G)$ with $2\leq |\E(S)|=k<d$. Consider an arbitrary subset $S\subset \FR_0$ with $|\E(S)|=k+1$, and denote $\G:=\G(S)$.
    
    If $|\V(S)|\geq |\E(S)|+2$, then the $\FR$-graph $\G$ is not connected, and hence we can write $\G$ as a nontrivial disjoint union $\G=\G_1\sqcup\G_2$. Note that $\M(\G)=\M(\G_1)\cup \M(\G_2)$ and $\M(\G_1)$ commutes with $\M(\G_2)$ by Lemma \ref{conjugate}(3). Then by Lemma \ref{basic lemma of FR-graph}, we have
    $$|\E(\M(\G_i))|\leq |\E(\M(\G))|-1= |\E(S)|-1\leq k.$$ 
    By the induction hypothesis, we obtain that $\M(\G_i)$ has a fixed point for $i=1,2$. Therefore, $\M(\G)$ has a fixed point by Corollary \ref{product maintain fix}.

    If $|\V(S)|\leq |\E(S)|+1\leq d+1$, then by Eq. (\ref{eq.FVE}), we have
     \begin{equation}\label{eq. n conjugate}
    n:=\lfloor\frac{2d+2-|\V(S)|}{\Cf(S)}\rfloor+1\geq \lfloor\frac{2d+1-|\E(S)|}{\E(S)}\rfloor+1\geq \lfloor\frac{2d+1}{|\E(S)|}\rfloor\geq 2.
     \end{equation}
     Since $\ell_i\geq 2d+2$ for all $i=1,\ldots, m$, we can construct $n$ conjugates of $\M(\G)$ as follows: first take $n$ permutations $\phi_1,\ldots,\phi_n\in \Perm(G)$, each fixing $\V(S)\setminus \CF(S)$ pointwise, such that $$\phi_j(\Cf(S))\cap \phi_{j'}(\Cf(S))=\emptyset, \quad \text{for all} ~1\leq j\neq j'\leq n;$$
     and then take $M_j:=\phi_j \M(\G)\phi_j^{-1}$. By Lemma \ref{conjugate}(1)(3), all the $n$ conjugates $M_j$ of $\M(\G)$ commute, i.e., $[s_j,s_{j'}]=1$ for all $s_j \in M_j$ and $s_{j'} \in M_{j'} (j \neq j')$.
    Moreover, by Lemma \ref{ad function2}, we have 
    $$\lfloor\frac{2d+1}{|\E(S)|} \rfloor(|\E(S)|-1) >d,$$ 
    where $3\leq |\E(S)|\leq d$. By the induction hypothesis, every subset of $\M(\G)$ with $|\E(S)|-1$ elements has a fixed point. Combining the above inequality, it follows that any finite subset of $\M(\G)$ has a fixed point by the conjugate bootstrap (Corollary \ref{Conjugate Bootstrap}), and so does $\M(\G)$ as $\generate{\M(\G)}$ is finitely generated. 
\end{proof}

We now consider some special cases when $|\E(S)|= d+1$.

\begin{lem}\label{V leq d+1 and E=d+1} 
 Let $S$ be as in Proposition \ref{solvable graph lemma} with $|\E(S)|= d+1$. If one of the following holds: 
 \begin{enumerate}
     \item $|\V(S)|\leq d+1$;
     \item $|\V(S)|=d+2$,  $\G(S)$ is a tree, and there exists a vertex in $\G(S)$ serving as the common origin of at least two directed edges,
 \end{enumerate}
 then $\M(\G(S))$ has a fixed point in $X$.
\end{lem}

\begin{proof}
    Note that in case (1), we have $|\Cf(S)|\leq d+1$ by Eq. (\ref{eq.FVE}), while in case (2), $|\Cf(S)|\leq d$ by hypothesis. Consequently, in both cases, $$\lfloor\frac{2d+2-|\V(S)|}{\Cf(S)}\rfloor+1\geq 2.$$ 
    Then by a similar discussion following Eq. (\ref{eq. n conjugate}), there exist at least $2$ commuting conjugates of $\M(\G(S))$. 
   Note that $2d>d$, and by Lemma \ref{|E(S)|leq d}, every $d$-element subset in $\M(\G(S))$ has a fixed point. Therefore, by the conjugate bootstrap (Corollary \ref{Conjugate Bootstrap}), $\M(\G(S))$ has a fixed point in $X$.
\end{proof}

We can immediately get the following corollary on $T_v^\G$ (see Definition \ref{Rooted trees}).

\begin{cor}\label{T have fixed pt}
Let $S$ be as in Proposition \ref{solvable graph lemma} with $|\E(S)|\leq d+1$. Then for any vertex $v$ in $\G:=\G(S)$ satisfying $v\notin \Post^\G(v)$,  $\M(T_v^\G)$ has a fixed point in $X$.
\end{cor}

\begin{proof}
    Now there are at most $d+1$ edges in $T_v^\G$. Then by Lemma \ref{|E(S)|leq d} and Lemma \ref{V leq d+1 and E=d+1}(2), $\M(T_v^\G)$ has a fixed point in $X$.
\end{proof}

\begin{lem}\label{G is a tree}
    Let $S$ be as in Proposition \ref{solvable graph lemma} with $|\E(S)|\leq d+1$. If $\G:=\G(S)$ is a tree, then $\M(\G)$ has a fixed point in $X$.
\end{lem}
\begin{proof}
If $|\E(S)|\leq d$, then $\M(\G)$ has a fixed point in $X$ by Lemma \ref{|E(S)|leq d}. We now assume $|\E(S)|= d+1$ in the following. Since $\G$ is a tree, there exists a vertex $v$ that is the origin of some edges but not the terminus of any edge. It follows that $v\notin \Post^\G(v)$. Furthermore, if at least two edges emanate from $v$, then $\M(\G)$ has a fixed point in $X$ by Lemma \ref{V leq d+1 and E=d+1}(2). We therefore assume in what follows that exactly one edge $e_0$ emanates from $v$.
    
    We first claim that for any edge $e'\in \Ecal(\G)\setminus \{e_0\}$, the edge group $G_{e'}=(G_{o(e')}, G_{\tau(e')})$ normalizes $\generate{\M(T_v^\G)}$, i.e., for any $\phi=(G_{o(e')}, c) \in G_{e'}$ with $c \in G_{\tau(e')}$ and any edge $e$ in $T_v^\G$, we have
    \begin{equation}\label{eq. Ge normalise Ge'}
        \phi^{-1}G_e\phi \subset \generate{\M(T_v^\G)}.
    \end{equation}
    Indeed, observe that the origin of $e$ in $T_v^\G$ is $v$, and $e_0$ is the unique edge emanating from $v$ in $\G$. We therefore have $v\neq o(e'), \tau(e')$, and
    $$G_e=(G_v,G_{\tau(e)}), \quad G_v \neq G_{o(e')}, ~G_{\tau(e')}.$$ 
    Then there are two cases:
    \begin{itemize}
        \item if $G_{o(e')}\neq G_{\tau(e)}$, then $G_{e'}$ and $G_{e}$ commute by Lemma \ref{conjugate}(3) and hence Eq. (\ref{eq. Ge normalise Ge'}) holds;
        \item  if $G_{o(e')}=G_{\tau(e)}$, then for any $\phi \in G_{e'}$, by Lemma \ref{conjugate}(4),
        $$\phi^{-1} G_e\phi \subset \generate{G_e\cup G_{e''}} \subset \generate{\M(T_v^\G)},$$
        where $e''$ is the edge in $T_v^\G$ emanating from $v$ to $\tau(e')$. Hence Eq. (\ref{eq. Ge normalise Ge'}) holds as well.
    \end{itemize}
    
    Now let $\G'$ be the subgraph of $\G$ obtained by deleting the edge $e_0$. Then
    both sets $\M(\G')$ and $\M(T^\G_v)$ have a fixed point in $X$ by Lemma \ref{|E(S)|leq d} and Corollary \ref{T have fixed pt} respectively. Combining the above claim with Proposition \ref{Normalize maintain fixed point}, 
    $\M(\G')\cup \M(T^\G_v)$ has a fixed point in $X$. Therefore, $\M(\G) \subset \M(\G')\cup \M(T^\G_v)$ has a fixed point in $X$.
\end{proof}

We have now arrived at our main goal.

\begin{proof}[\textbf{Proof of Proposition \ref{solvable graph lemma}}]
 For any subset $S\subset \FR_0(G)$, the $\FR$-graph $\G:=\G(S)$ is finite, with all its edge groups finitely generated. It follows that $\generate{\M(\G)}$ is finitely generated. Thus, it suffices to prove that any finite subset $S'\subset \M(\G)$ has a fixed point in $X$. Furthermore, by Lemma \ref{basic lemma of FR-graph}(1), $S'\subset \M(\G')$ where $\G'=\G(S')$ is the corresponding $\FR$-graph generated by $S'$. Thus, to prove that $S'$ has a fixed point in $X$, it suffices to prove that $\M(\G')$ has a fixed point in $X$. Moreover, Corollary \ref{Helly type theorem} implies that we only need to consider the case when $|S'|\leq d+1$. Now $\G'$ has at most $d+1$ edges. The case when $|\E(S')|\leq d$ has been proved in Lemma \ref{|E(S)|leq d}. To prove the case when $|\E(S')|=d+1$, we discuss in the following three subcases: 
   
   (1) $\G'$ is not connected. Then each component of $\G'$ has at most $d$ edges and the edge groups in distinct components commute. By Corollary \ref{product maintain fix} and Lemma \ref{|E(S)|leq d}, $\M(\G')$ has a fixed point in $X$.

   (2) $\G'$ is connected but not a tree. Then $|\V(S')|\leq d+1$ and hence $\M(\G')$ has a fixed point in $X$ by Lemma \ref{V leq d+1 and E=d+1}(1).

   (3) $\G'$ is a tree. Then by Lemma \ref{G is a tree}, $\M(\G')$ has a fixed point in $X$. This completes the proof. 
\end{proof}

\section{\texorpdfstring{$\phi$-torsion and Property (III)}{phi-torsion and Property (III)}}\label{sect 6}

Since the definition of Property (I) is relative to the action on the space, we introduce the following Property (III) only relying on groups.
\begin{defn}[$\phi$-torsion and Property (III)]
    Let $\phi\in \Aut(G)$ be an automorphism of a group $G$. An element $g\in G$ is said to be $\phi$-\emph{torsion} if there exists a positive integer $m$ such that 
    \begin{equation}\label{def. phi-order}
       \phi(g)\phi^2(g)\cdots\phi^{m}(g)=e. 
    \end{equation}
   The smallest such positive integer $m$ is called the \emph{$\phi$-order} of $g$. Moreover, we say $G$ satisfies \emph{Property (III)} if there exists an automorphism $\phi\in\Aut(G)$ and a generating set $\s_G$ such that all elements in $\s_G$ are $\phi$-torsion.
\end{defn}

A group $G$ is called \emph{invertible} if it has a generating set $\s_G$ such that the map $\tau(s) = s^{-1}$ for $s\in \s_G$ induces an involution $\tau\in \Aut(G)$, i.e., $\tau^2 = \id$.

Note that if $G$ is finitely generated, then it admits a finite generating set $\s_G'\subset \s_G$ for any generating set $\s_G$. Without loss of generality, we assume that $\s_G$ is finite when $G$ is finitely generated in the above definitions.
If, in addition, $G$ satisfies Property (III) with $\phi\in\Aut(G)$, then $\phi$ must be of finite order. Indeed, if $s\in \s_G$, applying $\phi^{-1}$ to Eq. (\ref{def. phi-order}), we have 
$s\phi(s)\ldots\phi^{m-1}(s)=e.$
Thus $s=\phi^{m}(s)$ and 
\begin{equation}\label{phi^M=id}
    \phi^M=\id,
\end{equation}
where $M$ is a common multiple of the $\phi$-orders of all generators in the finite generating set $\s_G$.

\begin{prop}\label{PIII groups} 
The following hold.
\begin{enumerate}
    \item  Finite groups satisfy Property (III).
    \item  Invertible groups, including closed orientable surface groups, Artin groups and Coxeter groups, all satisfy Property (III).
    \item  The solvable Baumslag-Solitar groups $BS_{1,n}(n\geq 2)$ do not satisfy Property (III).
    \item  Property (III) is preserved under direct and free products but not under semidirect products.
\end{enumerate}
\end{prop}
\begin{proof}
    (1) Any finite group $G$ satisfies Property (III) with $\phi=\id$ and $\s_G=G$.

    (2) Given an invertible group $G$ with $\s_G$ and the involution $\tau$, we have $$\tau(s)\tau^2(s) = s^{-1}s = e$$ for every $s \in \s_G$ and hence $G$ satisfies Property (III). 
    
    Let $G$ be a closed orientable surface group, i.e., $G=\pi_1(\Sigma_n)$ where $\Sigma_n$ is a closed orientable surface with genus $n\geq 2$. Then $G$ admits the following symmetric presentation
    $$\pi_1(\Sigma_n)=\generate{c_1,\ldots,c_{2n} \mid c_1\cdots c_{2n}c_1^{-1}\cdots c_{2n}^{-1}},$$
    and the map $\tau(c_i)=c_i^{-1}$ induces an involution. We obtain that $G$ is also invertible (see \cite{WZZ} for more information on symmetric presentations of surface groups). Artin groups (e.g. free groups, abelian groups) and Coxeter groups are invertible obviously.

    (3) Note that $$BS_{1,n}=\generate{a,t \mid tat^{-1}=a^n}=\langle\!\langle a \rangle\!\rangle\rtimes \langle t\rangle\cong \Z[1/n]\rtimes \Z.$$ Since $a=t^{-1}a^nt=(t^{-1}at)^n$, we can denote $t^{-1}at$ by $a^{1/n}$ and any element in $BS(1,n)$ can be written as $a^vt^u$ where $v\in \Z[1/n]$ and $u\in \Z$.  Furthermore, by \cite[Lemma 3.1]{ON}, any automorphism $\varphi:BS_{1,n} \to BS_{1,n}$ must be of the form $$\varphi(a)=t^{-k}a^lt^k=a^{l/n^k}  \text{ and } ~\varphi(t)=t^{-p}a^qt^{p+1}=a^{q/n^p}t$$ with integers $k,p\geq 0$ and $l \neq 0$ dividing some power of $n$.
   
    For any finite generating set $\s$ of $BS(1,n)$, there must exist an element $s=a^{v_0}t^{u_0}\in \s$ with $u_0\neq 0$. Then for any positive integer $m$, $$\phi(s)\cdots \phi^m(s)=w(a)t^{mu_0}\neq e,$$
    where $w(a)$ lies in the normal closure $\langle\!\langle a \rangle\!\rangle$ of $a$. This gives the result.

    (4) Let $G_i (i=1,2)$ be two groups satisfying Property (III), namely, there exist automorphisms $\phi_i\in \Aut(G_i)$ and generating sets $\s_i$ for $G_i$ such that all elements in $\s_i$ are $\phi_i$-torsion. 
    
    For the direct product $G_1\times G_2$, take its generating set as $\s_\times:=\{(s_1, e), (e, s_2)\mid s_i\in \s_i\}$. Then $\Phi:=\phi_1\times \phi_2$ 
    is an automorphism of $G_1\times G_2$ such that all elements in $\s_\times$ are $\Phi$-torsion.
    
    For the free product $G_1\ast G_2$, take its generating set as $\s_\ast:=\s_1\sqcup \s_2$. Then $\Psi:=\phi_1\ast \phi_2$ defined by $\Psi|_{\s_i}=\phi_i,$
    is an automorphism of $G_1\ast G_2$ such that all elements in $\s_\ast$ are $\Psi$-torsion.

    From the preceding discussions, direct and free products preserve Property (III). By contrast, item (3) shows that semidirect products do not.
\end{proof}

We introduce the following lemma to prove Theorem \ref{main2} in Section \ref{sect 7}.

\begin{lem}\label{property3 imply property1}
     Let $G=\ast_{i=1}^m \ast _{j=1}^{\ell_i}G_{i,j}$ be a free product as in Eq. (\ref{free decomposition}) and let $\Aut(G)$ act by isometries on a complete $\CAT0$ space $X$ with $\dim(X)\leq d$. If $\ell_i\geq 2d+2$ and $G_{i,j}$ satisfies Property (III) for each $i$, then
     \begin{itemize}
         \item $G$ satisfies Property (I);
         \item $H^1(\Aut(G),\Z)\cong \Pi_{i=1}^m H^1(\Aut(G_{i,1}),\Z)$.
     \end{itemize}
\end{lem}
\begin{proof}
    Let $A\neq B$ be free factors.
    Then $B$ has property (III), i.e., $B$ has a finite generating set $\s_B$ such that each element in $\s_B$ is $\phi$-torsion. 
    Since $\ell_i\geq 2d+2$, we can choose $d+1$ commuting conjugates of $(A,B)\rtimes \Aut(B)$ and obtain an embedding  
    \begin{equation}\label{eq. Lemma property3 imply property1}
        ((A,B)\rtimes \Aut(B))\wr S_{d+1}\hookrightarrow \Aut(G)=\FR(G)\rtimes \Pi_{i=1}^{m}(\Aut(G_{i,1})\wr S_{\ell_i})
    \end{equation}
     by Theorem \ref{Gilbert}. 
    By the finiteness of $\s_B$ and Eq. (\ref{phi^M=id}), we can take $M$ as a common multiple of the $\phi$-orders of all generators in $\s_B$ such that $\phi^M=\id \in \Aut(B)$. Then for any $b\in \s_{B}$, we have
        \begin{align}
            (\phi(A,b))^{M}&=\phi(A,b)\phi^{-1}\phi^{2}(A,b)\phi^{-2}\cdots\phi^{M}(A,b)\phi^{-M}\phi^{M}  \label{Eq. Property (III)}\\
            &=(A,\phi(b))(A,\phi^2(b))\cdots (A,\phi^{M}(b)) \notag \\ 
            &=(A,\phi(b)\phi^2(b)\cdots\phi^{M}(b))\notag \\ 
            &=\id\in (A,B)\rtimes \Aut(B), \notag 
        \end{align}
    where the second ``='' follows from Lemma \ref{conjugate}(1).
    Thus both $\phi(A,b)$ and $\phi$ are of finite order and by Lemma \ref{close under conjugation, radical, multiplication}, they lie in $ \Elp_X((A,B)\rtimes \Aut(B))$, and hence $(A,b)\in \generate{\phi(A,b),\phi}$ has a fixed point and so does $(A,B)=\generate{\{(A,b)\mid b\in \s_{B}\}}$. Therefore, the first statement holds.

    Since the short exact sequence
    $$1\to \FR(G)\to \Aut(G)\to Q:=\Pi_{i=1}^{m}(\Aut(G_{i,1})\wr S_{\ell_i})\to 1$$
    splits, we obtain the following decomposition of the first cohomology group:
    $$H^1(\Aut(G),\Z)\cong H^1(\FR(G),\Z)^{Q}\times H^{1}(Q,\Z),$$ where $$H^1(\FR(G),\Z)^{Q}=\{\chi \in {\rm Hom}(\FR(G), \Z)\mid \phi\cdot \chi = \chi ~\text{for all}~\phi \in Q\}.$$ 
    It suffices to prove $H^1(\FR(G),\Z)^Q=0$ as $H^{1}(Q,\Z)\cong \Pi_{i=1}^mH^{1}(\Aut(G_{i,1}),\Z)$.
    Indeed, given $[\chi]\in H^1(\FR(G),\Z)^{Q}$, i.e., $$\chi: \FR(G) \to \Z\quad  {\rm and} \quad \phi\cdot \chi(x)=\chi(\phi x\phi^{-1})=\chi(x) $$
    for any $\phi\in Q$ and $x\in \FR(G)$,
    we have
    $$0=\chi((\phi(A,b))^{M})\overset{Eq.~(\ref{Eq. Property (III)})}=\chi(\phi(A,b)\phi^{-1})+\chi(\phi^{2}(A,b)\phi^{-2})+\cdots+\chi(\phi^{M}(A,b)\phi^{-M})=M\chi(A,b),$$
    for any $b\in \s_B$.
    Thus, $\chi(A,b)=0$, which finishes the proof.
\end{proof}

\section{Fixed-point dimensions and some applications}\label{sect 7}

We first introduce the following definition from Bridson \cite{Bridson(Hellytype)}.

\begin{defn}
     The \emph{fixed-point dimension} of a group $G$, denoted by $\Fixdim(G)$, is the largest integer $d$ such that $G$ has a fixed point, whenever $G$ acts by isometries on a complete $\CAT0$ space $X$ with $\dim(X)\leq d$.
\end{defn}

We denote $\Fixdim(G)=\infty$ if $\Fixdim(G)\geq d$ for every integer $d$. For instance, the fixed-point dimension of any finite group is infinite.

\begin{proof}[\textbf{Proof of Corollary \ref{main thm 2}}]
    Let $\Aut(G)$ act by isometries on a complete $\CAT0$ space $X$ with $\dim(X)\leq d$. 
    Since $\Fixdim(G_{i,j})\geq d$ and $\Fixdim(\Aut(G_{i,j}))\geq d$ imply Properties (I) and (II) respectively, the conclusion follows from Theorem \ref{main thm}.
\end{proof}

\begin{proof}[\textbf{Proof of Theorem \ref{main2}}]
    The conclusion follows by combining Lemma \ref{property3 imply property1}, Proposition \ref{property 2 equal finite abelianization} and Theorem \ref{main thm}.
\end{proof}

As an application of Theorem \ref{main2}, we have
\begin{cor}
    Let 
    $G=\overbrace{\Z^{n_1}\ast\cdots\ast\Z^{n_1}}^{\ell_1}\ast
    \cdots \ast \overbrace{\Z^{n_m}\ast\cdots\ast\Z^{n_m}}^{\ell_m}~ (n_i\neq n_j \geq 2 \textit{ for } i\neq j)$. Then $$\Fixdim(\Aut(G))\geq \min_{1\leq i\leq m}\{\lfloor\frac{\ell_i}{2}\rfloor-1\}.$$
    Namely, any isometric action of $\Aut(G)$ on a complete $\CAT0$ space $X$ with $\dim(X)\leq d$ has a fixed point if each $\ell_i\geq 2d+2$.
\end{cor}
\begin{proof}
    By Proposition \ref{PIII groups} and Lemma \ref{property3 imply property1}, $\Z^{n_i}$ satisfies Property (III) and $$H^1(\Aut(G),\Z)=\Pi_{i=1}^m H^1(\Aut(\Z^{n_i}),\Z)=0$$ respectively. Note that $\Aut(\Z^{n_i})$ is finitely generated.
    Thus, the conclusion follows from Theorem \ref{main2}.
\end{proof}

For a closed orientable surface $\Sigma_{n}$ with genus $n\geq 1$, Bridson \cite{Bridson(MCG)} proved that
any semi-simple action of its mapping class group $\Mod(\Sigma_n)$ on any complete $\CAT0$ space of dimension less than $n$ has a fixed point. Our Corollary \ref{surface fixed point} considers automorphism groups of free products of surface groups acting by isometries on complete $\CAT0$ spaces.

\begin{cor}[Corollary \ref{surface fixed point}]
    Let $$G=\overbrace{\pi_1(\Sigma_{n_1})\ast\cdots\ast\pi_1(\Sigma_{n_1})}^{\ell_1}\ast
    \cdots \ast \overbrace{\pi_1(\Sigma_{n_m})\ast\cdots\ast \pi_1(\Sigma_{n_m})}^{\ell_m}\quad (n_i\neq n_j\geq 1 {~for~} i\neq j).$$
    Then any isometric action of $\Aut(G)$ on a complete $\CAT0$ space $X$ with $\dim(X)\leq d$ has a fixed point if each $\ell_i\geq 2d+2$.
\end{cor}
\begin{proof}
We first claim that $$H^1(\Aut(\pi_1(\Sigma_{n_i})),\Z)=0$$ for every $n_i\geq 1$. 
Indeed, if $n_i=1$, $\Aut(\pi_1(\Sigma_1))\cong \mathrm{GL}(2,\Z)$, then its first cohomology is trivial. If $n_i\geq 2$, then by the Dehn-Nielsen-Baer theorem (cf. \cite[Section 8.2.7]{Farb(MCG)}),
\begin{equation}\label{Aut=Mod pm}
    \Aut(\pi_1(\Sigma_{n_i})) \cong \mathrm{Mod}^{\pm}(\Sigma_{n_i,1})
\end{equation}
is finitely generated,
where $\Sigma_{n_i,1}$ denotes the surface with a marked point. Note that the orientation-preserving mapping class group $\mathrm{Mod}(\Sigma_{n_i,1})$ is an index-two subgroup of $\mathrm{Mod}^{\pm}(\Sigma_{n_i,1})$. By the Birman exact sequence (cf. \cite[Section 4.2.1]{Farb(MCG)}),
$$
1\longrightarrow \pi_1(\Sigma_{n_i})
\longrightarrow \mathrm{Mod}(\Sigma_{n_i,1})
\longrightarrow \mathrm{Mod}(\Sigma_{n_i})
\longrightarrow 1,
$$
the five-term exact sequence in group homology gives
$$
H_2(\mathrm{Mod}(\Sigma_{n_i}),\mathbb Z)
\longrightarrow
H_1(\pi_1(\Sigma_{n_i}),\mathbb Z)_{\mathrm{Mod}(\Sigma_{n_i})}
\longrightarrow
H_1(\mathrm{Mod}(\Sigma_{n_i,1}),\mathbb Z)\longrightarrow
H_1(\mathrm{Mod}(\Sigma_{n_i}),\mathbb Z)
\longrightarrow 0.
$$
However, the symplectic action of $\mathrm{Mod}(\Sigma_{n_i})$ on $H_1(\pi_1(\Sigma_{n_i}),\mathbb Z)$ implies that 
$$H_1(\pi_1(\Sigma_{n_i}),\mathbb Z)_{\mathrm{Mod}(\Sigma_{n_i})}=0.$$ 
Thus, we have 
$$H_1(\mathrm{Mod}(\Sigma_{n_i,1}),\mathbb Z)\cong
H_1(\mathrm{Mod}(\Sigma_{n_i}),\mathbb Z)$$ and they are both finite (cf. \cite[Sections 5.1.2 and 5.1.3]{Farb(MCG)}).
It follows that $H^1(\Aut(\pi_1(\Sigma_{n_i})),\Z)=0$ by Eq. (\ref{Aut=Mod pm}), i.e., the claim holds.

By Proposition \ref{PIII groups}, $\pi_1(\Sigma_{n_i})$ satisfies Property (III). Furthermore, by Lemma \ref{property3 imply property1}, we have $$H^1(\Aut(G),\Z) \cong \Pi_{i=1}^m H^1(\Aut(\pi_1(\Sigma_{n_i})),\Z)=0.$$
Then $\Aut(G)$ has a fixed point in $X$ by Theorem \ref{main2}.
\end{proof}

\appendix
\section{General cases}\label{section general}

In this appendix, we consider the group $G\ast F_r$, where $G=\ast_{i=1}^m\ast_{j=1}^{\ell_i} G_{i,j}$ is defined in Eq. (\ref{free decomposition}) and $F_r$ is a free group of rank $r\geq 1$, and show the following.
\begin{thm}\label{main thm 3}
    Suppose that $\Aut(G\ast F_r)$ acts by isometries on a complete $\CAT0$ space $X$ with $\dim (X)\leq d$, and the following hold:
    \begin{itemize}
        \item $\ell_i=0$, or $\ell_i\geq 2d+2$;
        \item $r=0$, or $r\geq 3s$ if $d=2s-1$, or $r\geq 3s+2$ if $d=2s$, for some integer $s$;
        \item $G_{i,j}$ satisfies Property (III);
        \item $\Aut(G)$ is finitely generated and $H^{1}(\Aut(G),\Z)=0$.
    \end{itemize}
    Then $\Aut(G\ast F_r)$ has a fixed point in $X$.
\end{thm}

Note that Theorem \ref{main thm 3} recovers Theorem \ref{Bridson theorem} (\cite[Theorem A]{Bridson(Hellytype)}) when all $\ell_i=0$. When $r=0$, it coincides with Theorem \ref{main2}.

\subsection{\texorpdfstring{Automorphism group $\Aut(G\ast F_r)$}{Automorphism group Aut(G*Fr)}}
Let $\s_G$ be defined as in Section \ref{sect 2.1} and let
\begin{equation}\label{S_r and S_Gr}
   \quad \s_r:=\{a_1,\ldots,a_r\} \quad\text{ and }\quad \s_{G\ast F_r}:=\s_G\sqcup \s_r
\end{equation}
be generating sets of $F_r$ and $G\ast F_r$ respectively.
Given $a_k\in \s_r$, $s\in \s_G$, there are three types of automorphisms in $\Aut(G\ast F_r)$ as follows.
\begin{itemize}
    \item \emph{Right transvections}:
    $$\rho(a_k,s)(x):=\begin{cases}
    a_ks\quad & x=a_k\\
    x\quad & x\in \s_{G\ast F_r}\setminus \{a_k\}
    \end{cases};$$

    \item \emph{Left transvections}:
    $$\lambda(a_k,s)(x):=\begin{cases}
    sa_k\quad & x=a_k\\
    x\quad & x\in \s_{G\ast F_r}\setminus \{a_k\}
    \end{cases};$$
    
    \item \emph{Dehn twists}:
    $$\FR(G_{i,j},a_k)(x):=\begin{cases}
    a_k^{-1}xa_k\quad & x\in \s_{G_{i,j}}\\
    x \quad & x\in \s_{G\ast F_r}\setminus \s_{G_{i,j}}
    \end{cases}.$$
\end{itemize}
Similar to partial conjugations in Eq. (\ref{PC}), left and right transvections correspond to certain subgroups of $\Aut(G\ast F_r)$,
$$\rho(a_k,G_{i,j}):=\{\rho(a_k,g)\mid g \in G_{i,j}\}\cong G_{i,j},$$
$$\lambda(a_k,G_{i,j}):=\{\lambda(a_k,g)\mid g \in G_{i,j}\}\cong G_{i,j}.$$
We can also define the following subsets of $\Aut(G\ast F_r)$:
\begin{itemize}
    \item $\rho(F_r,G_{i,j}):=\{\rho(a_k,g)\mid \text{$a_k\in \s_r$, $g\in G_{i,j}$}\}$;
    \item $\lambda(F_r,G_{i,j}):=\{\lambda(a_k,g)\mid \text{$a_k\in \s_r$, $g\in G_{i,j}$}\}$;
    \item $\FR(G,F_r):=\{\FR(G_{i,j},a_k)\mid \text{$a_k\in \s_r$ and $G_{i,j}$ are free factors in Eq. (\ref{free decomposition})}\}$.
\end{itemize}

Take 
$$\s_{\SAut(F_r)}:=\{\lambda_{i,j}, ~\rho_{i,j}\mid \textit{$1\leq i\neq j\leq r$}\}\subset\Aut(G\ast F_r),$$
where 
$$\rho_{i,j}(a_k):=\begin{cases}
    a_ka_j\quad & k=i\\
    a_k\quad & k\neq i
    \end{cases},\quad
   \lambda_{i,j}(a_k):=\begin{cases}
    a_ja_k\quad & k=i\\
    a_k\quad & k\neq i
    \end{cases}$$ are Nielsen transformations and we denote    $$\Niel:=\{\rho_{1,2},\lambda_{1,2},~\ldots,~\rho_{r-1,r},\lambda_{r-1,r},~\rho_{r,1},\lambda_{r,1}\}.$$
    Let $\epsilon_1$ be an automorphism of $F_r$ that maps $a_1$ to $a_1^{-1}$ with other generators invariant. It is well known that $$\generate{\s_{\SAut(F_r)}}=\generate{\Niel}, \quad \Aut(F_r)=\generate{\s_{\SAut(F_r)}}\rtimes\generate{\epsilon_1}.$$
Then we get a normal subgroup of $\Aut(G\ast F_r)$ analogous to $\generate{\FR_0(G)}$ defined in Eq. (\ref{eq. FR0}):
$$\generate{\FR_0(G\ast F_r)}\triangleleft \Aut(G\ast F_r),$$
where the generating set
    \begin{align}
        &\FR_0(G\ast F_r):=\FR_0(G)\cup\s_{\SAut(F_r)}\cup\FR(G,F_r)\label{decomposition of FR_0(G*Fr)} \cup \bigcup_{i,j} (\rho(F_r,G_{i,j})\cup\lambda(F_r,G_{i,j})) \\ 
        \supset &\FR_1(G\ast F_r):=\FR_0(G)\cup \Niel 
        \cup \bigcup_{i,j}(\rho(a_1,G_{i,j})\cup \lambda(a_1,G_{i,j})\cup\FR(G_{i,j},a_1)).\notag
    \end{align}
Furthermore, we have the following technical fact by straightforward calculation.
\begin{lem}\label{decomposition of FR_1(G*Fr)}
    $\generate{\FR_1(G\ast F_r)}\subset\generate{\FR_0(G\ast F_r)}\subset\generate{\FR_0(G\ast F_r)}\rtimes\generate{\epsilon_1}=\generate{\FR_1(G\ast F_r)}\rtimes\generate{\epsilon_1}$.
\end{lem}

Similar to Theorem \ref{Gilbert}, we have that
$\generate{\FR_0(G\ast F_r)}$ is normalized by $\generate{\epsilon_1,\Fact(G),\Perm(G)}$ and that
\begin{equation}\label{Gilbert2}
    \Aut(G\ast F_r)=\generate{\FR_0(G\ast F_r)} \generate{\epsilon_1,\Fact(G),\Perm(G)}
\end{equation}
is finitely generated. For more information on automorphism groups of free products, see \cite{Fuchs}.

\subsection{\texorpdfstring{Extended $\FR$-graph $\G(S)$}{Extended FR-graph G(S)}}

Now we define a new directed graph with respect to every subset $S\subset\FR_0(G\ast F_r)$, which is a natural generalization of $\FR$-graphs in Section \ref{sect 4}. Recall that {\bf we replace $S$ with $S\setminus\id$ if $\id \in S$} in the following.

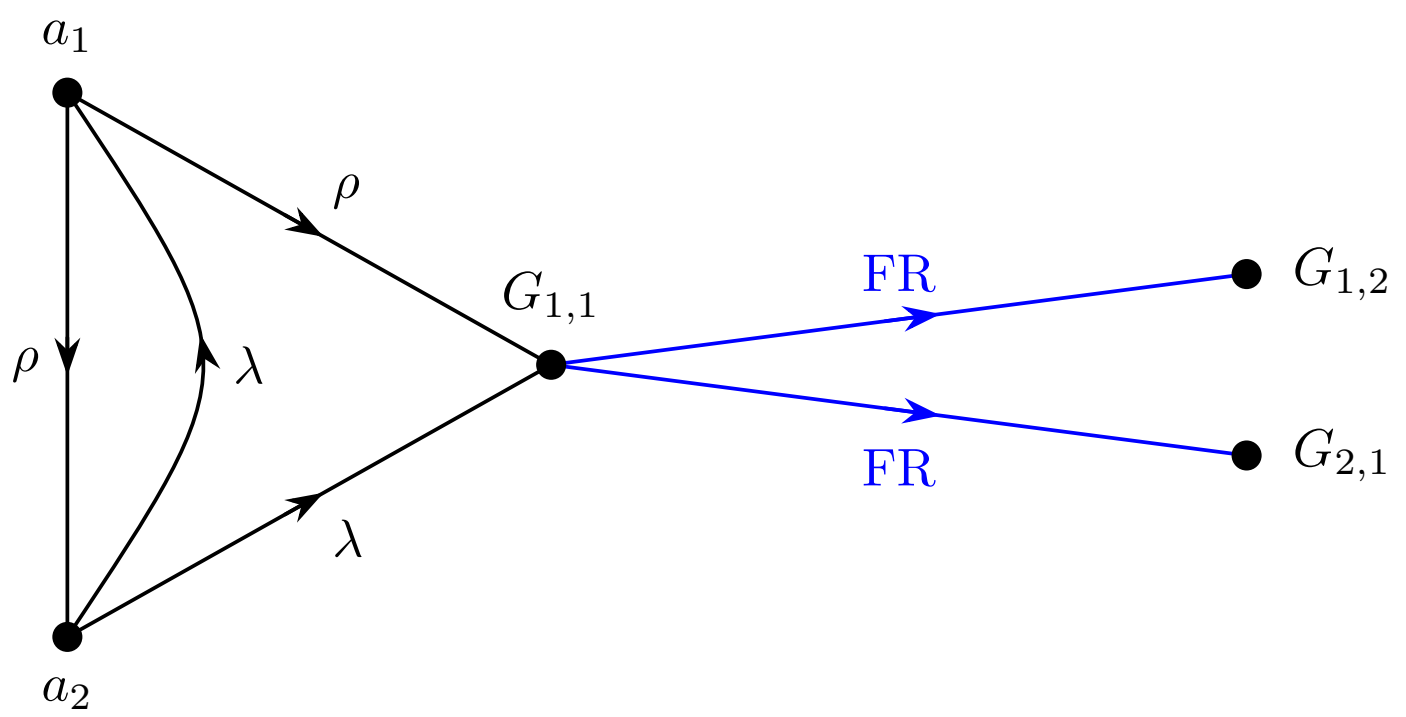
\begin{figure}[htbp]
\centering
\begin{tikzpicture}[
    line width=0.7pt,
    line cap=round,
    line join=round,
    vertex/.style={circle, fill=black, inner sep=2pt},
    edge/.style={
        -{Stealth[length=2.5mm,width=1.7mm]},
        line width=0.7pt
    }
]

\tikzset{
  midarrow/.style={
    postaction={
      decorate,
      decoration={
        markings,
        mark=at position 0.55 with {\arrow{Stealth[length=2.5mm,width=1.7mm]}}
      }
    }
  }
}

% Vertices
\coordinate (a1)   at (0,3.6);
\coordinate (a2)   at (0,0);
\coordinate (G11)   at (3.2,1.8);
\coordinate (G12) at (7.8,2.4);
\coordinate (G21) at (7.8,1.2);

% Edges
\draw (a2) -- node[midway, left=4pt, fill=white, inner sep=1pt]
    {$\rho$}(a1);
\draw[edge]
    ($(a1)!0.45!(a2)$) -- ($(a1)!0.52!(a2)$);

\draw[midarrow] 
  (a2) .. controls (1.2,1.8) and (1.2,1.8) .. (a1);
\node at (1.2,1.8) {$\lambda$};

\draw (a1) -- node[midway, above=7pt, right=3pt, fill=white, inner sep=1pt]
    {$\rho$}(G11);
\draw[edge]
    ($(a1)!0.45!(G11)$) -- ($(a1)!0.53!(G11)$);

\draw (a2) -- node[midway, below=7pt, right=3pt, fill=white, inner sep=1pt]
    {$\lambda$}(G11);
\draw[edge]
    ($(a2)!0.45!(G11)$) -- ($(a2)!0.53!(G11)$);

\draw[blue] (G11) -- node[midway, above=4pt, fill=white, inner sep=1pt]
    {\blue $\FR$}(G12);
\draw[blue, edge]
    ($(G11)!0.48!(G12)$) -- ($(G11)!0.56!(G12)$);

\draw[blue] (G11) -- node[midway, below=6pt, fill=white, inner sep=1pt]
    {\blue $\FR$}(G21);
\draw[blue, edge]
    ($(G11)!0.48!(G21)$) -- ($(G11)!0.56!(G21)$);

% Solid vertices
\foreach \v in {a1,a2,G11,G12,G21} {
    \node[vertex] at (\v) {};
}

% Labels
\node[above=4pt] at (a1) {$a_1$};
\node[below=4pt] at (a2) {$a_2$};
\node[above=4pt] at (G11) {$G_{1,1}$};

\node[right=5pt] at (G12) {$G_{1,2}$};
\node[right=5pt] at (G21) {$G_{2,1}$};

\end{tikzpicture}
\caption{An extended $\FR$-graph $\G(S)$ for 
$S =\{\rho_{1,2},\lambda_{2,1},\rho(a_1,h)\}\cup\lambda(a_2,G_{1,1})\cup(G_{1,1},G_{1,2})\cup(G_{1,1},G_{2,1})$, where $h\in G_{1,1}$.}
    \label{fig:extended FR-graph}
\end{figure}

\begin{defn}[Extended $\FR$-graph]
    For any subset $S\subset \FR_0(G\ast F_r)$, we define a finite directed graph $\bar\G(S)$ labeled by $S$ (see Figure \ref{fig:extended FR-graph} for an example), called the \emph{extended $\FR$-graph generated by} $S$, as follows.
    \begin{itemize}
        \item Vertex set $\Vcal(\bar\G(S))$: vertices are divided into two types,
        
        (v1): each $a_i\in \s_r$ denotes a vertex;
        
        (v2): each free factor $G_{i,j}$ denotes a vertex.

        \item Edge set $\Ecal(\bar\G(S))$: edges are divided into four types and each type of edges corresponds to some elements contained in $S$. More precisely, 

        (e1): if $\rho_{i,j}\in S$ (resp. $\lambda_{i,j} \in S$), then draw a directed edge from $a_i$ to $a_j$ with label $\rho$ (resp. $\lambda$);
        
        (e2): if $\rho(a_k,G_{i,j})\cap S\neq \emptyset$ (resp. $\lambda(a_k,G_{i,j})\cap S\neq \emptyset$), then draw a directed edge from $a_k$ to $G_{i,j}$ with label $\rho$ (resp. $\lambda$);
        
        (e3): if $\FR(G_{i,j},a_k)\in S$, then draw a directed edge from $G_{i,j}$ to $a_k$ with label $\FR$;
        
        (e4): if $(G_{i,j}, G_{i',j'})\cap S\neq \emptyset$, then draw a directed edge from $G_{i,j}$ to $G_{i',j'}$ with label $\FR$.
    \end{itemize}  
    Given two vertices $u,v$, there is at most one edge from $u$ to $v$ with each label. We write $\G(S)$ after omitting all isolated vertices of $\bar \G(S)$ and denote $\G=\G(S)$ when no confusion arises. 
\end{defn}

%--------------------------------------------------------------------- -------------------

Note that if $S\subset \FR_0(G)$, then the extended $\FR$-graph $\G(S)$ is exactly the $\FR$-graph with each edge labeled by $\FR$. Thus, we will not distinguish $\FR$-graphs and extended $\FR$-graphs in the following.

Similar to $\FR$-graphs, we can define vertex groups, edge groups for extended $\FR$-graphs and adopt their notations $G_v,G_e$ as in Section \ref{sect 4}. In fact, any vertex group is a free factor $G_{i,j}$ or $\generate{a_k}$ ($a_k\in \s_r$) and any edge group is one of the following forms: $$\generate{\rho_{k,l}},~ \generate{\lambda_{k,l}}; \quad \rho(a_k,G_{i,j}),~ \lambda(a_k,G_{i,j}); \quad \FR(G_{i,j}, \generate {a_k}); \quad (G_{i,j},G_{i',j'}).$$
We also denote the union of all edge groups as $\M(\G(S))$. 

By a similar argument as in Lemma \ref{conjugate}(3), we have the following facts on extended $\FR$-graphs.
\begin{lem}\label{appendix commute}
    Let $S_1,S_2\subset \FR_0(G\ast F_r)$.
    \begin{enumerate} 
        \item If $|S_1|=1$, then $\G(S_1)$ is a segment. 
        \item If $\G(S_1)\cap\G(S_2)=\emptyset$, or $\G(S_1)\cap\G(S_2)$ consists of some vertices which do not emit edges in $\Ecal(\G(S_1\cup S_2))$, then $S_1$ and $S_2$ commute pairwise.
    \end{enumerate}
\end{lem}

\begin{defn}[Directed triangle, post-vertex set, cutting point]\label{directed triangle, post-vertex set, cutting point}
    Let $\G$ be an extended $\FR$-graph. 
    \begin{enumerate}
        \item A \emph{directed triangle} $\Delta(A,B,C)\subset \G$ consists of three points $A,B,C \in \Vcal(\G)$ and three directed edges $e_{A,B}, e_{B,C}, e_{A,C}$ such that $e_{A,B}$ and $e_{A,C}$ share the same label. The three edge groups $G_{e_{A,B}}, G_{e_{A,C}}, G_{e_{B,C}}$ satisfy the following \emph{triangle conjugacy relation}: $$\phi^{-1}\generate{G_{e_{A,B}}}\phi \subset \generate{G_{e_{A,B}} \cup G_{e_{A,C}}},\quad\forall \phi\in G_{e_{B,C}}.$$
        
        \item Similar to the $\FR$-graph, for a vertex $v$ in $\G$, we can also define the \emph{post-vertex set}: $$\Post^\G(v):=\{u\in \Vcal(\G)\mid \text{there exists a directed path from $v$ to $u$ in $\G$}\}.$$
        
        \item A vertex $v$ in $\G$ is called a \emph{cutting point}, if $\G$ can be  decomposed along $v$ into two subgraphs $\G_1,\G_2$ such that 
        $$\G_1\cup \G_2=\G, \quad \G_1\cap \G_2=\{v\} \quad \text{and}\quad \Post^{\G_1}(v)=\emptyset.$$

        \item Given a cutting point $v\in \Vcal(\G)$, we define a new extended $\FR$-graph $\tilde\G_1$ for $\G_1$ as follows: 
            \begin{itemize}
                \item vertex set $\Vcal(\tilde\G_1):=\Vcal(\G_1)\cup \Post^{\G}(v)$;
                \item edge set $\Ecal(\tilde\G_1):=\Ecal(\G_1) \cup  \{e_{w,u} \mid e_{w,v} \in \Ecal(\G_1), u \in \Post^{\G}(v) \}$, where $e_{w,u}$ is a new directed edge from $w$ to $u$ sharing the same label as $e_{w,v}$.
            \end{itemize}
    \end{enumerate}
    For example, in Figure \ref{fig:extended FR-graph}, the vertices $a_1$, $a_2$ and $G_{1,1}$ form a directed triangle $\Delta(a_1,a_2,G_{1,1})$, and $$\Post^{\G}(a_2)=\{a_1, a_2, G_{1,1}, G_{1,2}, G_{2,1}\}, \quad \Post^{\G}(G_{1,1})=\{G_{1,2}, G_{2,1}\}.$$ The graph $\G(S)$ can be decomposed along the cutting point $G_{1,1}$ into two graphs $$\G_1=\Delta(a_1,a_2,G_{1,1}),\quad  \G_2=\G((G_{1,1},G_{1,2})\cup(G_{1,1},G_{2,1})).$$
    The graph $\tilde \G_1$ for $\G_1=\Delta(a_1,a_2,G_{1,1})$ is shown in Figure \ref{fig: tilde G_1}.
\end{defn}

\begin{figure}[htbp]
\begin{tikzpicture}[
    line width=0.7pt,
    line cap=round,
    line join=round,
    vertex/.style={circle, fill=black, inner sep=2pt},
    edge/.style={
        -{Stealth[length=2.5mm,width=1.7mm]},
        line width=0.7pt
    }
]

\tikzset{
  midarrow/.style={
    postaction={
      decorate,
      decoration={
        markings,
        mark=at position 0.55 with {\arrow{Stealth[length=2.5mm,width=1.7mm]}}
      }
    }
  }
}

% Vertices
\coordinate (a1)   at (6,3.6);
\coordinate (a2)   at (6,0);
\coordinate (G11)   at (9.2,1.8);
\coordinate (G12) at (3.9,1.8);
\coordinate (G21) at (1.8,1.8);

% Edges
\draw (a1) -- node[midway, right=4pt, fill=white, inner sep=1pt]
    {$\rho$}(a2);
\draw[edge]
    ($(a1)!0.45!(a2)$) -- ($(a1)!0.52!(a2)$);

\draw[midarrow] 
  (a2) .. controls (7.2,1.8) and (7.2,1.8) .. (a1);
\node at (7.2,1.8) {$\lambda$};

\draw (a1) -- node[midway, above=7pt, right=3pt, fill=white, inner sep=1pt]
    {$\rho$}(G11);
\draw[edge]
    ($(a1)!0.45!(G11)$) -- ($(a1)!0.53!(G11)$);

\draw (a2) -- node[midway, below=7pt, right=3pt, fill=white, inner sep=1pt]
    {$\lambda$}(G11);
\draw[edge]
    ($(a2)!0.45!(G11)$) -- ($(a2)!0.53!(G11)$);

\draw (a1) -- node[midway, below=4pt, fill=white, inner sep=1pt]
    {$\rho$}(G12);
\draw[edge]
    ($(a1)!0.48!(G12)$) -- ($(a1)!0.56!(G12)$);

\draw (a1) -- node[midway, left=4pt, above=4pt, fill=white, inner sep=1pt]
    {$\rho$}(G21);
\draw[edge]
    ($(a1)!0.48!(G21)$) -- ($(a1)!0.56!(G21)$);

\draw (a2) -- node[midway, above=4pt, fill=white, inner sep=1pt]
    {$\lambda$}(G12);
\draw[edge]
    ($(a2)!0.48!(G12)$) -- ($(a2)!0.56!(G12)$);

\draw (a2) -- node[midway, left=4pt, below=4pt, fill=white, inner sep=1pt]
    {$\lambda$}(G21);
\draw[edge]
    ($(a2)!0.48!(G21)$) -- ($(a2)!0.56!(G21)$);

% Solid vertices
\foreach \v in {a1,a2,G11,G12,G21} {
    \node[vertex] at (\v) {};
}

% Labels
\node[above=4pt] at (a1) {$a_1$};
\node[below=4pt] at (a2) {$a_2$};
\node[right=2pt] at (G11) {$G_{1,1}$};

\node[left=2pt] at (G12) {$G_{1,2}$};
\node[below=2pt] at (G21) {$G_{2,1}$};

\end{tikzpicture}
\caption{}
    \label{fig: tilde G_1}
\end{figure}
%----------------------------------------------------------------------------------------

We now give some properties of the union $\M(\G(S))$ of all edge groups.
\begin{lem}\label{triangle lemma}
If an extended $\FR$-graph $\G$ is decomposed into $\G_1$ and $\G_2$ along a cutting point $v$, then
\begin{enumerate}
    \item $\M(\G)=\M(\G_1)\cup \M(\G_2)\subset \M(\tilde \G_1)\cup \M(\G_2)$;
    \item $\M(\G_2)$ normalizes $\generate{\M(\tilde \G_1)}$.
\end{enumerate}
\end{lem}
\begin{proof}
    Assertion (1) is obvious. Now we prove assertion (2).
    
    For any $\psi\in \M(\tilde \G_1)$, $\phi \in \M(\G_2)$, let  $e_\psi \in \Ecal(\tilde \G_1)$, $e_\phi \in\Ecal(\G_2)$ denote the edges whose edge groups contain $\psi$, $\phi$ respectively. 
    Since $\Vcal(\tilde \G_1)\cap \Vcal(\G_2)\subset\{v\}\cup \Post^{\G}(v)$, we have $o(e_\psi)\neq o(e_{\phi}), \tau(e_\phi)$. Thus, there are only two cases.
    
        If the endpoints of $e_\psi$ and $e_\phi$ are disjoint or $\tau(e_\psi)=\tau(e_\phi)$, then $\psi$ and $\phi$ commute by Lemma \ref{appendix commute}(2), and hence $\phi^{-1} \psi \phi\in \M(\tilde \G_1)$.
        
        If $\tau(e_\psi)=o(e_\phi)$, then the triangle conjugacy relation implies $\phi^{-1} \psi \phi\in \generate{G_{e_\psi}\cup G_{\tilde e}}\subset \generate{\M(\tilde \G_1)}$, where $\tilde e\in \Ecal(\tilde \G_1)$ is an edge from $o(e_\psi)$ to $\tau(e_\phi)$ sharing the same label as $e_\psi$. 
\end{proof}

By Lemma \ref{triangle lemma},
we get the following useful result.
\begin{lem}\label{normalize criterion}
     Let $\Aut(G\ast F_r)$ act by isometries on a complete $\CAT0$ space $X$, and let $\G$ be an extended $\FR$-graph decomposed into $\G_1$ and $\G_2$ along a cutting point $v$.
     If $\M(\tilde{\G_1})$ and $\M(\G_2)$ both have a fixed point in $X$, then $\M(\G)$ has a fixed point in $X$.
\end{lem}

\subsection{Proof of Theorem \ref{main thm 3}}

We first consider the edge groups.
\begin{prop}\label{one edge extended FR-graph}
    Let $\Aut(G\ast F_r)$ act by isometries on a complete $\CAT0$ space $X$ with $\dim(X)\leq d$ and the following hold:
    \begin{itemize}
        \item $\ell_i=0$ or $\ell_i\geq 2d+2$;
        \item $r=0$, or $r\geq 3s$ if $d=2s-1$, or $r\geq 3s+2$ if $d=2s$, for some integer $s$;
        \item $G_{i,j}$ satisfies Property (III).
    \end{itemize}
    Then every edge group of $\G(\FR_0(G\ast F_r))$ has a fixed point in $X$.
\end{prop}
\begin{proof}
    An edge group has one of the following types: $$\generate{\rho_{i,j}},~\generate{\lambda_{i,j}},~ \rho(a_k,G_{i,j}),~\lambda(a_k,G_{i,j}), ~\FR(G_{i,j},\generate{a_k}), ~(G_{i,j},G_{i',j'}).$$ 
    Note that $$\generate{\rho_{i,j}},\generate{\lambda_{i,j}}\hookrightarrow\Aut(F_r)\hookrightarrow \Aut(G\ast F_r)$$ and thus $\generate{\rho_{i,j}},\generate{\lambda_{i,j}}$ both have a fixed point by Theorem \ref{Bridson theorem}. 
    
    Since $G_{i,j}$, $\generate{a_k}$ satisfy Property (III), the remaining cases are similar to the proof of the first statement in Lemma \ref{property3 imply property1} by replacing $(A,B)\rtimes \Aut(B)$ in Eq. (\ref{eq. Lemma property3 imply property1}). For example, $$(\rho(a_k,G_{i,j})\rtimes \Aut(G_{i,j}))\wr S_{d+1}\hookrightarrow \Aut(G\ast F_r)$$ acts by isometries on $X$. Then, by a similar argument below Eq. (\ref{eq. Lemma property3 imply property1}), $\rho(a_k,G_{i,j})$ has a fixed point as $G_{i,j}$ satisfies Property (III).
\end{proof}

Before a key proposition, we recall the Ample Duplication Criterion \cite[Definition 5.1, Theorem 5.2]{Bridson(Hellytype)}.

\begin{defn}[\cite{Bridson(Hellytype)}, Ample Duplication]\label{AD definition}
    Given a group $\Gamma$ and positive integers $d$ and $k_0$, we say that a finite generating set $\A$ for a subgroup $\Lambda<\Gamma$ has \emph{ample duplication for dimension $d$, with base $k_0$}, if there is an \emph{ample duplication function} $f:\mathbb{N}\to \mathbb{N}$, i.e., $f$ satisfies the following conditions:
    \begin{itemize}
        \item each subset $S\subset \A$ of cardinality $|S|>k_0$ can be written as a disjoint union $S=S_1\sqcup S_2$ where $S_i$ are nonempty and $\generate{S_1}$ normalizes $\generate{S_2}$, or else there are at least $f(|S|)$ commuting conjugates of $\generate{S}$ in $\Gamma$;
        \item $d<(k-1)f(k)$ for $k=k_0+1,\ldots,\min\{d+1,|\A|\}$.
    \end{itemize}
\end{defn}

\begin{thm}[\cite{Bridson(Hellytype)}, Ample Duplication Criterion]\label{ADC}
    Let $\Gamma$ be a group acting by isometries on a complete $\CAT0$ space $X$ with $\dim(X)\leq d$, and let $\Lambda<\Gamma$ be a subgroup with a finite generating set $\A$ that has ample duplication for dimension $d$, with base $k_0$. If $\generate{S}$ has a fixed point for every $S\subset \A$ with $|S|\leq k_0$, then $\Lambda$ has a fixed point in $X$.
\end{thm}
\begin{remark}\label{remark of ADC}
    In the Ample Duplication Criterion, it suffices to require the first condition in Definition \ref{AD definition} only for those subsets $S\subseteq A$ with $|S|>k_0$ that do not have a fixed point in $X$. Indeed, subsets that already have a fixed point can be omitted in the inductive argument on $|S|$.
\end{remark}

\begin{figure}[htbp]

\begin{tikzpicture}[scale=0.95, line width=0.9pt, >=Stealth]

\path[use as bounding box] (-6.3,-5.0) rectangle (5.3,0.8);

\begin{scope}[xshift=0cm]

% ---------- arrow in the middle ----------
\tikzset{
  midarrow/.style={
    postaction={
      decorate,
      decoration={
        markings,
        mark=at position 0.55 with {\arrow{Stealth[length=2.5mm,width=1.7mm]}}
      }
    }
  }
}

% ---------- helper command: lens-shaped double edge ----------
\newcommand{\lens}[2]{
  \draw[midarrow]
    (#1) .. controls ($(#1)+(0.40,0.30)$) and ($(#2)+(-0.40,0.30)$) .. (#2);
  \draw[midarrow]
    (#1) .. controls ($(#1)+(0.40,-0.30)$) and ($(#2)+(-0.40,-0.30)$) .. (#2);
}

% ---------- top row nodes ----------
\coordinate (L0) at (-6.0,0);
\coordinate (L1) at (-4.75,0);
\coordinate (L2) at (-3.50,0);

\coordinate (C0) at (-1.65,0);  % a_r
\coordinate (C1) at (0,0);      % a_1
\coordinate (C2) at (1.65,0);   % a_2

\coordinate (R1) at (3.50,0);   % a_\ell
\coordinate (R2) at (4.75,0);   % a_{\ell+1}

% ---------- draw top dots ----------
\foreach \P in {L0,L1,L2,C0,C1,C2,R1,R2}{
  \fill (\P) circle (2.2pt);
}

% ---------- top lens edges ----------
\lens{L0}{L1}
\node at (-5.5,0.44) {$\rho$};
\node at (-5.5,-0.44) {$\lambda$};
\lens{L1}{L2}
\node at (-4.25,0.44) {$\rho$};
\node at (-4.25,-0.44) {$\lambda$};

\node at (-2.58,0.02) {$\cdots$};

\lens{C0}{C1}
\node at (-1.0,0.44) {$\rho$};
\node at (-1.0,-0.44) {$\lambda$};
\lens{C1}{C2}
\node at (0.85,0.44) {$\rho$};
\node at (0.85,-0.44) {$\lambda$};

\node at (2.58,0.02) {$\cdots$};

\lens{R1}{R2}
\node at (4.0,0.44) {$\rho$};
\node at (4.0,-0.44) {$\lambda$};

% ---------- top labels ----------
\node[below=0.32cm] at (L0) {$a_t$};
\node[below=0.32cm] at (L1) {$a_{t+1}$};
\node[below=0.32cm] at (L2) {$a_{t+2}$};

\node[below=0.32cm] at (C0) {$a_r$};
\node[above=0.32cm] at (C1) {$a_1$};
\node[below=0.32cm] at (C2) {$a_2$};

\node[below=0.32cm] at (R1) {$a_\ell$};
\node[below=0.32cm] at (R2) {$a_{\ell+1}$};

% ---------- bottom nodes ----------
%
\coordinate (Gi) at (-1.45,-3.68);
\coordinate (Gj) at ( 1.45,-3.68);

% outer oval
\draw[dashed, very thick] (0,-4.02) ellipse [x radius=3.65, y radius=0.76];

\fill (Gi) circle (2.3pt);
\fill (Gj) circle (2.3pt);

\node[below=0.12cm] at (Gi) {$G_{i,j}$};
\node[below=0.12cm] at (Gj) {$G_{i',j'}$};
\node at (3.4,-3.3) {$\FR$-graph};

% =====================================================
% left bundle: two arrows downward, one arrow upward
% =====================================================

% outer downward arrow

\draw[midarrow]
  (C1) .. controls (-1.3,-1.00) and (-2.24,-2.38) .. (Gi);
\node at (-1.7,-1.55) {$\rho$};

% middle downward arrow
\draw[midarrow]
  (C1) .. controls (-0.68,-1.18) and (-1.52,-2.48) .. (Gi);
\node at (-1.15,-1.65) {$\lambda$};

% inner upward arrow
\draw[blue, very thick, midarrow]
  (Gi) .. controls (-0.80,-2.38) and (-0.28,-1.08) .. (C1);
\node at (-0.3,-1.90) {\blue \tiny$\FR$};

% =====================================================
% right bundle: two arrows downward, one arrow upward
% =====================================================

% outer downward arrow
\draw[midarrow]
  (C1) .. controls (1.3,-1.00) and (2.24,-2.38) .. (Gj);
\node at (1.7,-1.55) {$\rho$};

% middle downward arrow
\draw[midarrow]
  (C1) .. controls (0.68,-1.18) and (1.52,-2.48) .. (Gj);
\node at (1.15,-1.65) {$\lambda$};

% inner upward arrow
\draw[blue, very thick, midarrow]
  (Gj) .. controls (0.80,-2.38) and (0.28,-1.08) .. (C1);
\node at (0.3,-1.90) {\blue \tiny$\FR$};
\end{scope}

\end{tikzpicture}
\caption{A standard extended FR-graph with edges lying in $\FR_1$, whose vertices are divided into two layers: the (v1)-type vertices on the top layer and the (v2)-type vertices on the bottom layer. The bottom layer is precisely an ordinary FR-graph.}
    \label{fig:extended FR graph}

\end{figure}
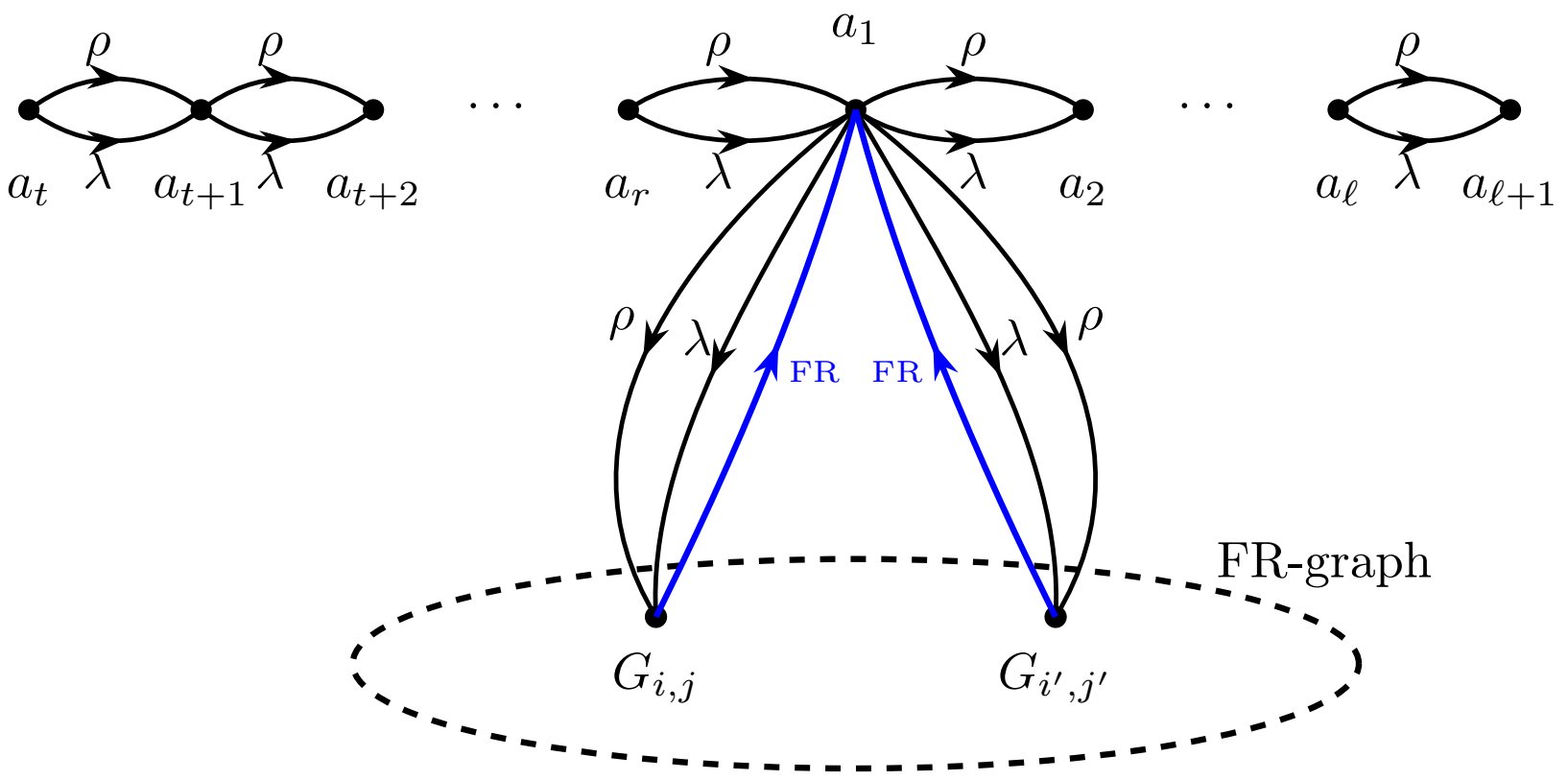

Let $\FR_0(G\ast F_r)$ and $\FR_1(G\ast F_r)$ be defined as in Eq. (\ref{decomposition of FR_0(G*Fr)}).

\begin{prop}\label{main prop in appendix}
    Under the hypothesis of Proposition \ref{one edge extended FR-graph}, $\FR_0(G\ast F_r)$ has a fixed point in $X$.
\end{prop}
\begin{proof}

    If all $\ell_i=0$ or $r=0$, then the result comes from the fact that $\Fixdim(\Aut(F_r))\geq d$ (Theorem \ref{Bridson theorem}) or $(\FR(G))$ has a fixed point by combining Proposition \ref{solvable graph lemma} and Lemma \ref{property3 imply property1}. Thus we assume each $\ell_i\neq 0$ and $r\neq 0$ in the following. 
    
    Note that $\generate{\epsilon_1}$ has a fixed point as it is a finite group. If we can prove $\FR_1(G\ast F_r)$ has a fixed point, then by Corollary \ref{product maintain fix}, $\generate{\FR_1(G\ast F_r)}\rtimes\generate{\epsilon_1}$ has a fixed point and hence by Lemma \ref{decomposition of FR_1(G*Fr)}, its subgroup $\generate{\FR_0(G\ast F_r)}$ has a fixed point. 

    To prove that, since $\generate{\FR_1(G\ast F_r)}$ is finitely generated, by Corollary \ref{Helly type theorem}, we only need to prove that any subset 
    \begin{equation}\label{S leq d+1}
        S\subset \FR_1(G\ast F_r)\quad \text{with}\quad |S|\leq d+1
    \end{equation}
    has a fixed point. 
    Furthermore, recall that $S\subset \M(\G(S))$, it suffices to prove $\M(\G(S))$ has a fixed point. The extended $\FR$-graph $\G=\G(S)$ is shown as in Figure \ref{fig:extended FR graph}.
    
    Without loss of generality, we assume that $\G$ is connected. Otherwise, $\G$ consisting of disjoint connected components $\G_i$ implies that $\M(\G_i)$ commute pairwise. Thus $\M(\G)=\cup\M(\G_i)$ has a fixed point by Corollary \ref{product maintain fix}, provided each $\M(\G_i)$ does.
    The proof proceeds by induction on the number of vertices from $\s_r$ in Eq. (\ref{S_r and S_Gr}), i.e., the number of
    $$ \Vcal_F(\G):=\{a_i\in \s_r\mid \text{$a_i$ is a (v1)-type vertex in $\G$}\} \subset \Vcal(\G).$$

    Case $|\Vcal_F(\G)|=0$ is trivial. If $|\Vcal_F(\G)|=1$, then the finite set $S$ has ample duplication for dimension $d$, with base $k_0=1$ in the sense of Remark \ref{remark of ADC}. Indeed, for any $S'\subset S$ of cardinality $|S| \geq |S'|> k_0=1$ with $\G(S')$ connected (else decomposed into connected components), if $$S'\subset \Aut(G)\hookrightarrow\Aut(G\ast F_r),$$
    then $S'$ has a fixed point by Proposition \ref{solvable graph lemma}. Otherwise, there are at least $m(S') : =\min\{r,\lfloor 2d+2/|S'|\rfloor\}$ commuting conjugates of $S'$. In addition,
    $$d< (|S'|-1)m(S'),\quad  k_0+1\leq |S'|\leq \min\{d+1,|S|\}.$$
    Note that $S'$ has a fixed point if $|S'|=1$ by Proposition \ref{one edge extended FR-graph}. Thus using the Ample Duplication Criterion (Theorem \ref{ADC}), $S$ has a fixed point. 
    Therefore, $\M(\G)$ has a fixed point by Corollary \ref{Helly type theorem} as $\generate{\M(\G)}$ is finitely generated.
    
    We assume the result holds for all $|\Vcal_F(\G)|=k\leq d+1$, and now consider the case 
    $$|\Vcal_F(\G)|=k+1\leq d+2.$$
    There are two remaining cases.
    
    \textbf{Case ($\Vcal(\G)=\Vcal_F(\G)$)}. Then $\M(\G)\subset \Aut(F_r)$ and it has a fixed point as in the case of all $\ell_i=0$ in the beginning of our proof.
    
    \textbf{Case ($\Vcal(\G)\supsetneq\Vcal_F(\G)$)}. Then by a simple calculation, we have
    $$ |\Vcal_F(\G)|< |\Vcal(\G)|\leq |\Ecal(\G)|+1\leq |S|+1\overset{Eq.~(\ref{S leq d+1})}{\leq} d+2 \leq r.$$
    Thus, there is some $a_i\in \s_r$ but $a_i \notin \Vcal_F(\G)$. If $a_1\notin \Vcal_F(\G)$, then by the definition of $\FR_1(G\ast F_r)$, the extended $\FR$-graph $\G$ is not connected (see Figure \ref{fig:extended FR graph}), contradicting the connectivity hypothesis. Therefore, we have $a_1\in \Vcal_F(\G)$ and two subcases by considering whether $a_r \in \s_r$ lies in $\Vcal_F(\G)$.
    
    \textbf{Subcase ($a_r\in \Vcal_F(\G)$)}. Then $\{\rho_{r,1},\lambda_{r,1}\}\cap S\neq \emptyset$. Let $t$ be the number such that $a_t\in \Vcal_F(\G)$ and $a_{t-1}\notin \Vcal_F(\G)$, see Figure \ref{fig:extended FR graph}.
    Then $a_{t+1}$(if $t=r$, taking $a_{t+1}=a_1$) is a cutting point and $\G$ can be decomposed into two parts along $a_{t+1}$:
    $$\G_1:=\G(\{\rho_{t,t+1},\lambda_{t,t+1}\}\cap S),\quad\G_2:=(\G \setminus \G_1) \cup \{a_{t+1}\}.$$ 
    Thus we obtain 
    $$\M(\tilde \G_1)\subset \generate{\bar\M_1},$$ 
    where the extended $\FR$-graph $\tilde \G_1$ is constructed from $\G_1$ (see Definition \ref{directed triangle, post-vertex set, cutting point}(4)) and
    \begin{align*}
        \bar\M_1:=&\{\rho_{t,k},\lambda_{t,k}\mid a_k\in \Post^{\G}(a_{t+1})\}\cup \{\rho_{t,t+1},\lambda_{t,t+1}\}
        \\
        &\cup \{\text{generators of $\rho(a_t,G_{i,j})$ and $\lambda(a_t,G_{i,j})$}\mid G_{i,j}\in \Post^{\G}(a_{t+1})\}.
    \end{align*}
    
    We claim that the finite set $\bar \M_1$ has ample duplication for dimension $d$, with base $k_0=1$. 
    Indeed, for any $S'\subset \bar \M_1$ of cardinality $|\bar \M_1|\geq |S'|>k_0=1$ with $\G(S')$ connected (else decomposed into connected components), there are at least $m(S'):=\max\{2,r-|S'|\}$ commuting conjugates of $S'$. In addition,
    $$d< (|S'|-1)m(S'),\quad k_0+1\leq |S'|\leq \min\{d+1,|\bar \M_1|\}.$$
    Note that $S'$ has a fixed point if $|S'|=1$ by Proposition \ref{one edge extended FR-graph}. Thus using the Ample Duplication Criterion (Theorem \ref{ADC}), $\bar \M_1$ has a fixed point and so does $\M(\tilde \G_1)$. Furthermore, $\M(\G_2)$ has a fixed point by induction. Then by Lemma \ref{normalize criterion}, $\M(\G)$ has a fixed point. 

    \textbf{Subcase ($a_r\notin \Vcal_F(\G)$)}. Since $|\Vcal_F(\G)|\geq 2$, we have that $a_{2}\in \Vcal(\G)$ is a cutting point and $\G$ can be decomposed into two parts along $a_{2}$:    
    $$\G_2:=\G(\{\rho_{i,i+1},\lambda_{i,i+1}\mid i=2,\ldots,r-2\}\cap S),\quad \G_1:=(\G \setminus \G_2)\cup\{a_{2}\}.$$
    Thus we obtain $\M(\tilde \G_1)\subset \generate{\bar\M_1}$, where
    $$\bar \M_1:= \{\rho_{1,i},\lambda_{1,i}\mid a_i\in \Vcal_F(\G),i\neq 1\}\cup  \{\text{generators of edge groups in $\G_1$}\}.$$
    Then by the same argument as above (now $S'\subset \Aut(G)$ or $m(S') : =\max\{2,\min\{r-|S'|,\lfloor2d+2/|S'|\rfloor\}\}$), $\M(\G)$ has a fixed point.
\end{proof}
\begin{proof}[\textbf{Proof of Theorem \ref{main thm 3}}]
    Since $\Aut(G)\hookrightarrow \Aut(G\ast F_r)$, we have that $\Fact(G)\subset \Aut(G)$ has a fixed point by Theorem \ref{main2}.
    Furthermore, since $\generate{\epsilon_1}$ and $\Perm(G)$ are finite groups, $$\generate{\epsilon_1,\Fact(G),\Perm(G)}=\generate{\epsilon_1}\times (\Fact(G)\rtimes \Perm(G))$$ has a fixed point by Corollary \ref{product maintain fix}. Therefore, by combining Proposition \ref{Normalize maintain fixed point}, Eq. (\ref{Gilbert2}) and Proposition \ref{main prop in appendix}, $\Aut(G\ast F_r)$ has a fixed point.
\end{proof}

\noindent\textbf{AI declaration.} Most of the manuscript was completed prior to the widespread adoption of generative AI tools. ChatGPT-5.6 was used to improve the editing and formatting of the manuscript, and to perform consistency checks throughout the text. The authors reviewed all AI-assisted suggestions incorporated into the manuscript and take full responsibility for its contents.

\end{document}